\documentclass[10pt,twoside, a4paper,reqno]{amsart}

\usepackage{amscd,amssymb,amsmath,graphicx,verbatim,mathrsfs,xcolor}
\usepackage[british]{babel}
\usepackage{wasysym}
\usepackage{hyperref}
\usepackage{enumitem}
\usepackage{microtype}
\usepackage[capitalize]{cleveref}
\usepackage{mathtools}
\usepackage{multicol}
\usepackage[margin=2.5cm]{geometry}
 
\usepackage[normalem]{ulem}
\usepackage[foot]{amsaddr}
\usepackage[only,llbracket,rrbracket]{stmaryrd}
\binoppenalty=\maxdimen
\relpenalty=\maxdimen

\usepackage{tikz,mathrsfs}
\usetikzlibrary{babel,arrows,decorations.markings,arrows.meta,decorations.pathmorphing,cd,patterns,backgrounds}
\tikzset{
	symbol/.style={
		draw=none,
		every to/.append style={
			edge node={node [sloped, allow upside down, auto=false]{$#1$}}}
	}}
\tikzset{>=stealth'}
\tikzcdset{arrow style=tikz, diagrams={>=stealth'}, row sep=6em, column sep=6em}
\def\arrowLengthDisplayStyle{4ex}
\def\arrowHeightDisplayStyle{.8ex}
\def\arrowSkipDisplayStyle{.5ex}
\def\arrowLengthTextStyle{3ex}
\def\arrowHeightTextStyle{.8ex}
\def\arrowSkipTextStyle{.4ex}
\def\arrowLengthScriptStyle{2.5ex}
\def\arrowHeightScriptStyle{.6ex}
\def\arrowSkipScriptStyle{.3ex}
\def\arrowLengthScriptScriptStyle{2ex}
\def\arrowHeightScriptScriptStyle{.4ex}
\def\arrowSkipScriptScriptStyle{.2ex}

\renewcommand{\to}{\arrow{->}}

\newcommand{\epi}{\arrow{->>}}
\newcommand{\embed}{\arrow{right hook->}}

\newcommand{\MakeTikzArrowWithSuperscriptSubscript}[4]
{
	\mathchoice
	{ 
		\hspace*{\arrowSkipDisplayStyle}
		\begin{tikzpicture}[baseline]
		\draw [#1] (0,\arrowHeightDisplayStyle) -- node [above] {$#2$} node [below] {$#3$} (#4 * \arrowLengthDisplayStyle, \arrowHeightDisplayStyle);
		\end{tikzpicture}
		\hspace*{\arrowSkipDisplayStyle}
	}
	{ 
		\hspace*{\arrowSkipTextStyle}
		\begin{tikzpicture}[baseline]
		\draw [#1] (0,\arrowHeightTextStyle) -- node [above] {$\scriptstyle #2$} node [below] {$\scriptstyle #3$} (#4 * \arrowLengthTextStyle, \arrowHeightTextStyle);
		\end{tikzpicture}
		\hspace*{\arrowSkipTextStyle}
	}
	{ 
		\hspace*{\arrowSkipScriptStyle}
		\begin{tikzpicture}[baseline]
		\draw [#1] (0,\arrowHeightScriptStyle) -- node [above] {$\scriptscriptstyle #2$} node [below] {$\scriptscriptstyle #3$} (#4 * \arrowLengthScriptStyle, \arrowHeightScriptStyle);
		\end{tikzpicture}
		\hspace*{\arrowSkipScriptStyle}
	}
	{ 
		\hspace*{\arrowSkipScriptScriptStyle}
		\begin{tikzpicture}[baseline]
		\draw [#1] (0,\arrowHeightScriptScriptStyle) -- node [above] {$\scriptscriptstyle #2$} node [below] {$\scriptscriptstyle #3$} (#4 * \arrowLengthScriptScriptStyle, \arrowHeightScriptScriptStyle);
		\end{tikzpicture}
		\hspace*{\arrowSkipScriptScriptStyle}
	}
}

\newcommand{\MakeTikzArrowWithCentralLabel}[3]
{
	\mathchoice
	{ 
		\hspace*{\arrowSkipDisplayStyle}
		\begin{tikzpicture}[baseline]
		\draw [#1] (0,\arrowHeightDisplayStyle) -- node [fill=white,inner sep=1pt] {$#2$} (#3 * \arrowLengthDisplayStyle, \arrowHeightDisplayStyle);
		\end{tikzpicture}
		\hspace*{\arrowSkipDisplayStyle}
	}
	{ 
		\hspace*{\arrowSkipTextStyle}
		\begin{tikzpicture}[baseline]
		\draw [#1] (0,\arrowHeightTextStyle) -- node [fill=white,inner sep=1pt] {$\scriptstyle #2$} (#3 * \arrowLengthTextStyle, \arrowHeightTextStyle);
		\end{tikzpicture}
		\hspace*{\arrowSkipTextStyle}
	}
	{ 
		\hspace*{\arrowSkipScriptStyle}
		\begin{tikzpicture}[baseline]
		\draw [#1] (0,\arrowHeightScriptStyle) -- node [fill=white,inner sep=1pt] {$\scriptscriptstyle #2$} (#3 * \arrowLengthScriptStyle, \arrowHeightScriptStyle);
		\end{tikzpicture}
		\hspace*{\arrowSkipScriptStyle}
	}
	{ 
		\hspace*{\arrowSkipScriptScriptStyle}
		\begin{tikzpicture}[baseline]
		\draw [#1] (0,\arrowHeightScriptScriptStyle) -- node [fill=white,inner sep=1pt] {$\scriptscriptstyle #2$} (#3 * \arrowLengthScriptScriptStyle, \arrowHeightScriptScriptStyle);
		\end{tikzpicture}
		\hspace*{\arrowSkipScriptScriptStyle}
	}
}

\def\arrow#1{\def\lastArrowStyle{#1}
	\futurelet\testchar\arrowMaybeStreched}
\def\arrowMaybeStreched{\ifx[\testchar \let\next\arrowStreched
	\else \let\next\arrowUnstreched \fi
	\next}

\def\arrowStreched[#1]{\def\lastArrowStrech{#1}
	\futurelet\testchar\arrowMaybeLabel}
\def\arrowUnstreched{\def\lastArrowStrech{1}
	\futurelet\testchar\arrowMaybeLabel}

\def\arrowMaybeLabel{\ifx^\testchar \let\next\arrowSuperscript
	\else \ifx_\testchar \let\next\arrowSubscript
	\else \ifx~\testchar \let\next\arrowCentralLabel
	\else \let\next\arrowNoLabel
	\fi
	\fi
	\fi
	\next}

\def\arrowSuperscript^#1{\def\lastArrowSuperscript{#1}
	\futurelet\testchar\arrowSuperMaybeSub}
\def\arrowSuperMaybeSub{\ifx_\testchar \let\next\arrowSuperscriptSubscript
	\else \let\next\arrowSuperscriptNoSubscript \fi
	\next}

\def\arrowSubscript_#1{\def\lastArrowSubscript{#1}
	\futurelet\testchar\arrowSubMaybeSuper}
\def\arrowSubMaybeSuper{\ifx^\testchar \let\next\arrowSubscriptSuperscript
	\else \let\next\arrowSubscriptNoSuperscript \fi
	\next}

\def\arrowSuperscriptSubscript_#1{\def\lastArrowSubscript{#1}
	\arrowDrawSupSub}
\def\arrowSuperscriptNoSubscript{\def\lastArrowSubscript{}
	\arrowDrawSupSub}
\def\arrowSubscriptSuperscript^#1{\def\lastArrowSuperscript{#1}
	\arrowDrawSupSub}
\def\arrowSubscriptNoSuperscript{\def\lastArrowSuperscript{}
	\arrowDrawSupSub}
\def\arrowNoLabel{\def\lastArrowSuperscript{}
	\def\lastArrowSubscript{}
	\arrowDrawSupSub}

\def\arrowCentralLabel~#1{\MakeTikzArrowWithCentralLabel{\lastArrowStyle}{#1}{\lastArrowStrech}}
\def\arrowDrawSupSub{\MakeTikzArrowWithSuperscriptSubscript{\lastArrowStyle}{\lastArrowSuperscript}{\lastArrowSubscript}{\lastArrowStrech}}

\tikzcdset{arrow style=tikz, diagrams={>=stealth'}, row sep=4em, column sep=4em}

\counterwithin*{equation}{subsubsection}

\makeatletter

\renewcommand\subsubsection{\@secnumfont}{\bfseries}%
\renewcommand\subsubsection{\@startsection{subsubsection}{3}
  \z@{.5\linespacing\@plus.7\linespacing}{-.5em}%
  {\normalfont\bfseries}}
  
  \makeatother
  
  \makeatletter
\newcommand{\crefnames}[3]{%
  \@for\next:=#1\do{%
    \expandafter\crefname\expandafter{\next}{#2}{#3}%
  }%
}
\makeatother

\newcommand{\nref}[1]{\nameref{#1} \ref{#1}}

\crefnames{part,chapter,section}{\S\!}{\S\S\!}

\newtheorem{proposition-definition}[theorem]{Proposition-Definition}

\theoremstyle{definition}

\makeatletter
\newtheorem*{rep@theorem}{\rep@title}
\newcommand{\newreptheorem}[2]{%
	\newenvironment{rep#1}[1]{%
		\def\rep@title{#2 \ref{##1}}%
		\begin{rep@theorem}}%
		{\end{rep@theorem}}}
\makeatother

\newreptheorem{theorem}{Theorem}
\newreptheorem{lemma}{Lemma}
\newreptheorem{corollary}{Corollary}
\newreptheorem{proposition}{Proposition}
\newreptheorem{setup}{Setup}

\setlist{leftmargin=15pt,labelindent=15pt}
\setlist[enumerate]{wide=\parindent, leftmargin=15pt, labelwidth=15pt, align=left,label=(\roman*)}

\counterwithin*{footnote}{subsection}

\mathchardef\mhyphen="2D 
\newlist{subenumerate}{enumerate}{1}
\setlist[subenumerate,1]{label=(\arabic*)}

\let\originalleft\left
\let\originalright\right
\renewcommand{\left}{\mathopen{}\mathclose\bgroup\originalleft}
\renewcommand{\right}{\aftergroup\egroup\originalright}

\newcommand\defeq{\stackrel{\mathrm{\mbox{\scriptsize{def}}}}{\,=\,}}
\newcommand{\sth}{\,\,\vert\,\,}

\newcommand{\Z}{\mathbb{Z}}

\DeclareMathOperator{\add}{add}
\DeclareMathOperator{\thick}{thick}

\newcommand{\Der}{\mathscr{D}}

\DeclareMathOperator{\per}{per}

\DeclareMathOperator{\twoper}{per_{[0,1]}}
\DeclareMathOperator{\Mod}{Mod}
\DeclareMathOperator{\fgMod}{\mathrm{mod}}

\renewcommand{\k}{\mathit{k}}

\newcommand{\A}{\mathscr{A}}
\newcommand{\B}{\mathscr{B}}
\renewcommand{\L}{\mathscr{L}}

\newcommand{\E}{\mathbb{E}}

\renewcommand{\i}{\mathscr{I}}
\newcommand{\p}{\mathscr{P}}

\newcommand{\q}{\mathscr{Q}}
\newcommand{\calq}{\mathcal{Q}}
\renewcommand{\r}{\mathscr{R}}
\newcommand{\calr}{\mathcal{R}}
\newcommand{\s}{\mathscr{S}}
\newcommand{\z}{\mathcal{Z}}
\newcommand{\x}{\mathcal{X}}

\newcommand{\exactstr}{\mathcal{S}}

\renewcommand{\c}{\mathcal{C}}

\newcommand{\tcmc}[1]{\mathfrak{W}(#1)}
\newcommand{\picspace}[1]{\mathfrak{X}\left( #1 \right)}
\newcommand{\picgroup}[1]{{G}\left(#1\right)}

\renewcommand{\c}{\mathcal{C}}

\newcommand{\rigid}{\mathrm{rigid}}
\newcommand{\silt}{\mathrm{silt}}
\DeclareMathOperator{\K}{K}

\DeclareMathOperator{\dgcat}{\bf dgcat}
\DeclareMathOperator{\Hqe}{\bf Hqe}
\newcommand{\zeroAus}{\bf 0\mhyphen\! Aus}

\DeclareMathOperator{\Ch}{Ch}
\newcommand{\op}{^{\mathrm{op}}}

\newcommand{\Ddgb}{\mathscr{D}_{\mathrm{dg}}^{\mathrm{b}}}
\newcommand{\classspace}{\mathrm{B}}
\newcommand{\bfg}{\mathbf{g}}

\DeclareMathOperator{\Fun}{Fun}

\newcommand{\proj}{\mathrm{proj}}
\newcommand{\inj}{\mathrm{inj}}
\newcommand{\projinj}{\mathrm{proj\mhyphen inj}}

\title[Silting reduction and picture categories]{Silting reduction and picture categories of 0-Auslander extriangulated categories}

\author[Erlend D. Børve]{Erlend D. Børve}
\keywords{Exact dg category, 0-Auslander extriangulated category, Bongartz completion, generalized concentric twin cotorsion pair, silting reduction, $\tau$-cluster morphism category, picture category, picture space, picture group}
\subjclass[2020]{16B50, 16E45, 18G35, 18G80}

\address{Institut Fourier, Université Grenoble Alpes, 100 Rue des mathématiques, 38610 Gières, France (Current address: Abteilung Mathematik, Department Mathematik/Informatik der Universität zu Köln, Weyertal 86--90, 50931 Cologne, Germany).}
\email{erlend.d.borve@gmail.com}
\thanks{The author acknowledges support from the French ANR grant CHARMS (ANR-19-CE40-0017-02).}
\begin{document}
\begin{abstract}
	Let $\mathcal{C}$ be an extriangulated category and let $\mathcal{R}\subseteq \mathcal{C}$ be a rigid subcategory. Generalizing Iyama--Yang silting reduction, we devise a technical condition $\textbf{(gCP)}$ on $\mathcal{R}$ which is sufficient for the Verdier quotient $\mathcal{C}/\mathrm{thick}(\mathcal{R})$ to be equivalent to an ideal quotient. In particular, the Verdier quotient $\mathcal{C}/\mathrm{thick}(\mathcal{R})$ will admit an extriangulation in such a way that the localization functor $L_{\calr}\colon \mathcal{C} \rightarrow \mathcal{C}/\mathrm{thick}(\mathcal{R})$ is extriangulated.
	When $\mathcal{C}$ is 0-Auslander, the condition $\textbf{(gCP)}$ holds for all rigid subcategories $\mathcal{R}$ admitting Bongartz completions. Furthermore, we prove that the Verdier quotient $\mathcal{C}/\mathrm{thick}(\mathcal{R})$ then remains 0-Auslander. 
	As an application, we define the picture category of a $0$-Auslander exact dg category $\mathscr{A}$ with Bongartz completions, which generalizes the notion of $\tau$-cluster morphism category. We show that the picture category of $\A$ is a cubical category, in the sense of Igusa. The picture group of $\A$ is defined as the fundamental group of its picture category. When $H_0\A$ is $\mathbf{g}$-finite, the picture group of $\A$ is finitely presented.
\end{abstract}
	
	\maketitle
	
	\tableofcontents
	
\newpage	
	
\section{Introduction}

\subsection{} Additive categorifications of cluster algebras \cite{BMRRT06,Ami09} have had a revolutionary impact on the representation theory of associative algebras. In classical tilting theory \cite{BB06,HR82}, the mutation of a tilting module at an indecomposable direct summand is not always defined, but by modelling this mutation procedure on the mutation of seeds in cluster algebras, one obtains mutationally complete generalizations of tilting theory. Notable frameworks include cluster tilting theory \cite{BMR07}, $\tau$-tilting theory \cite{AIR14} and two-term silting \cite{IJY14}. There is also a rich interaction with torsion theories, t-strucutres, co-t-strucutres, as well as other category-theoretic notions \cite{KY12,BY14}.

\subsection{}
The mutational completeness of support $\tau$-tilting pairs and two-term silting complexes gives rise to reduction techniques, known as $\tau$-tilting reduction \cite{Jas15, BM18w} and silting reduction \cite{IY18,IY20}, respectively.
One can think of the $\tau$-cluster morphism category as a categorification these reduction techniques.
The first definitions $\tau$-cluster morphism categories were confined to hereditary algebras of type $A$, where the definition is combinatorial in nature \cite{Igu14,CN17}. Later, the definition took a representation-theoretic form, extending the setting to representation-finite herediary algebras \cite{IT17}, then to $\tau$-tilting finite algebras \cite{BM18w}, and ultimately to all finite-dimensional algebras \cite{BH21}. When going beyond the hereditary case, the composition rule was defined in terms of support $\tau$-tilting reduction \cite[§3]{BM18w}, a class of bijections defined in a case-by-case manner. When composing maps, the number of cases multiplies, rendering the proof of the associative property rather involved \cite[§§5--10]{BM18w}. 

\subsection{}\label{frameworks}
There are two generalized approaches to $\tau$-cluster morphism categories that considerably simplify the proof of the associative property. In previous work, the author defined the $\tau$-cluster morphism category of a connective differential graded algebra with finite-dimensional homology in all degrees  \cite{Bor21}. In this homological approach, one defines the composition rule using Iyama--Yang silting reduction \cite{IY18,IY20}, whose functoriality leads to a short proof of associativity. The second framework is combinatorial, and is due to Schroll--Tattar--Treffinger--Williams \cite{STTW23} and Kaipel \cite{Kai23}. The former set of authors show that the $\tau$-cluster morphism category can be recovered from the wall-and-chamber structure, whereas Kaipel generalizes their work by developing the of notion partitioned fan and then extracting a category from the underlying geometry.

\subsection{Summary of results}
In \Cref{sec:exactdg}, we give a brief introduction to the theory of exact dg categories, following Chen \cite{Che24thesis,Che23thesis}. Very roughly speaking, one arrives at the definition by accommodating Keller's minimal axioms \cite[Appendix A]{Kel90} of Quillen exact categories \cite{Qui73} in the framework of differential graded categories. We recall in \Cref{Che23.6.19} how the zeroth homology category of an exact dg category admits the structure of an extriangulated category, in the sense of Nakaoka--Palu \cite{NP19}.
We give a brief overview of localization theory for extriangulated categories in \Cref{subsec:loc}. In \Cref{subsec:0Ausdef}, we recall the notions of rigid and silting subcategories of an extriangulated category, with emphasis on the case of 0-Auslander extriangulated categories, in the sense of Gorsky--Nakaoka--Palu \cite{GNP23}.

Our original contributions are in \Cref{subsec:siltred,subsec:tcmc}. In \Cref{subsec:siltred}, we generalize Iyama--Yang silting reduction \cite{IY18,IY20} to extriangulated categories, extending the scope provided by the current literature \cite{LZZZ21,PZ24,LZ22}. The following technical condition will be imposed in \Cref{ass:cotor}.
\subsubsection*{Condition \textbf{(gCP)}}\label{ass:cotorin}
Let $(\c,\E)$ be an extriangulated category. Consider the following condition on a rigid subcategory $\calr$ of $\c$: 
	\begin{enumerate}
		\item[]\label{ass:cotor0}  The ordered tuple $\big((\calr^{\vee}, \calr^{\perp}), ({^{\perp}\calr}, \calr^{\wedge})\big)$ is a generalized concentric twin cotorsion pair (henceforth shortened to gCTCP) in $(\c,\E)$, where we use the notation introduced in \nref{setup:siltred}.
	\end{enumerate}
\newpage	
Our main results will be the following.
\begin{reptheorem}{prop:0Aus_cotor}
Let $(\c,\E)$ be a 0-Auslander extriangulated category and let $\calr\subseteq \c$ be a rigid subcategory.
Then $\calr$ admits Bongartz completions (see \Cref{Bongartz}) precisely when $\big((\calr, \calr^{\perp}),({^{\perp}\calr}, \calr)\big)$ is a gCTCP in $(\c,\E)$. In particular, \nameref{ass:cotorin} holds for all rigid subcategories of $\c$ admitting Bongartz completions.
\end{reptheorem}
\begin{reptheorem}{cor:siltred0Aus} 
Let $(\c,\E)$ be an extriangulated category and let $\calr\subseteq \c$ be a rigid subcategory.
\begin{enumerate}
	\item\label{prop:siltredin} Suppose that $(\c,\E)$ has enough projectives or injectives, and assume that \nameref{ass:cotorin} holds for $\calr$.
Then the composite functor
	\begin{equation*}
	\begin{tikzcd}
		{\z_{\calr}}  \defeq \calr^{\perp} \cap {^{\perp}\calr} \arrow[r,hook] & \c \arrow[r,"L_{\calr}"] & \c/\thick(\calr),
	\end{tikzcd}
	\end{equation*}
	where $L_{\calr}$ is the Verdier localization (see \Cref{NOSloc}), induces an equivalence $\overline{\c}_{\calr} \defeq {\z_{\calr}}  / [\calr] \to \c/\thick(\calr)$. Consequently, the category $\c/\thick(\calr)$ can be equipped with an extriangulation\footnote{In general, a Verdier localization of an extriangualted category is only weakly extriangualted, cf. \Cref{NOSloc}.} such that $L_{\calr}$ carries the structure of extriangulated functor.
	\item\label{cor:siltred0Aus1in}  If $(\c,\E)$ and $\calr$ are as in \ref{prop:siltred}, then Verdier localization functor $ \c \to^{L_{\calr}} \c/\thick(\calr)$ induces a bijection
\begin{equation*}
\rigid_{\calr}(\c) \to \rigid(\c/\thick(\calr))
\end{equation*}
from the set of rigid subcategories of $\c$ containing $\calr$ to the set of rigid subcategories of $\c/\thick(\calr)$. 
It restricts to an isomorphism of posets
\begin{equation*}
\silt_{\calr}(\c) \to \silt(\c/\thick(\calr))
\end{equation*}
from the poset of silting subcategories of $\c$ containing $\calr$ to the poset of silting subcategories of $\c/\thick(\calr)$.\footnote{We refer to \Cref{AT22.5.12} for the definition of the partial orders in question.}
	\item\label{cor:siltred0Aus2in} If $(\c,\E)$ is 0-Auslander and admits Bongartz completions, then the statements in \ref{prop:siltredin} and \ref{cor:siltred0Aus1in} above hold (by \nref{prop:0Aus_cotor}), and the Verdier quotient $\c/\thick(\calr)$ is 0-Auslander and admits Bongartz completions.  	
\end{enumerate}
\end{reptheorem}

We spend \Cref{subsec:tcmc} defining the picture category $\tcmc{\A}$ of a 0-Auslander exact dg category $\A$ (with cofibrant mapping complexes) in which all rigid subcategories of the extriangulated category $H_0\A$ admit Bongartz completions. By considering the special case $\A=\per_{[0,1]}(A)$, namely the two-term dg category of a dg category $A$ with finite-dimensional homology in all degrees (two-term dg categories of finite-dimensional algebras are examples of such), one recovers the definition proposed in previous work \cite{Bor21}. The \textit{picture space} $\picspace{\A}$ of $\A$ will be defined as the classifying space of $\tcmc{\A}$, and the \textit{picture group} $\picgroup{\A}$ is the fundamental group of the picture space. If the extriangulated category $H_0\A$ is Krull--Schmidt and contains a silting object. then the picture category $\tcmc{\A}$ is a cubical category in the sense of Igusa, which is shown in \Cref{prop:HI18.2.14}. In \nref{prop:picgroup_pres}, we give a presentation of the picture group in terms of generators and relations, and show that $\picgroup{\A}$ is finitely presented if $\A$ is $\bfg$-finite.

\subsection{Notation and conventions}\label{notcon}
A Grothendieck universe $\mathfrak{U}$ is fixed throughout, and $\mathfrak{U}$-small sets will simply be called \textit{small}. We say that a category $\c$ is \textit{small} the set $\bigsqcup_{x,y\in\c} \c(x,y)$ is small. All sets and categories we fix in our results are small, unless declared as \textit{large}.
A small field $\k$ is fixed throughout.

We employ the homological indexing convention for complexes.
Let $\Mod(\k)$ denote the monoidal category of $\k$-modules, and $\Ch(\k)$ that of chain complexes thereof. 
A \textit{$k$-category} is an enriched category over $\Mod(\k)$, whereas a \textit{dg $k$-category} is enriched over $\Ch(\k)$. The differential in a dg category will be denoted $\partial$. 
By considering the underlying $\k$-categories, one defines a notion of smallness for (dg) $\k$-categories.

Given a dg $\k$-category $\A$, one can produce two $\k$-categories, denoted $Z_0\A$ and $H_0\A$, having the same objects as $\A$, with morphism spaces given by
\begin{align*}
&Z_0\A(-,-) \defeq \ker \big(\A(-,-)_0 \to^{\partial_0}  \A(-,-)_{-1}\big) \\ 
\text{and}\quad &H_0\A(-,-) \defeq Z_0\A(-,-)/ \mathrm{im}\big( \A(-,-)_1 \to^{\partial_1}  \A(-,-)_{0} \big),
\end{align*}
respectively. The $\k$-category $H_0\A$ will be called the \textit{zeroth homology $\k$-category} of $\A$. {We have that $H_0\A$ is small precisely when $\A$ is small.}
Since a dg subcategory of a dg $\k$-category $\A$ is uniquely determined by the induced subcategory of $H_0\A$, we will often implicitly descend to $H_0\A$ when defining dg subcategories of $\A$.

A dg functor $F\colon \A\to \B$ is called a \textit{quasi-equivalence} if the chain map
\begin{equation*}
F_{x,y}\colon \A(x,y) \to  \B(Fx,Fy)
\end{equation*}
is a quasi-isomorphism for all $x,y\in\A$ and the functor $H_0F\colon H_0\A \to H_0\B$ is an equivalence of $\k$-categories. 
The category of small dg $\k$-categories admits a cofibrantly generated model structure in which the weak equivalences are the quasi-equivalences \cite{Tab05}. This model category will be denoted $\dgcat_{\k}$, and its homotopy category will be denoted $\Hqe_{\k}$. 

\subsection{Acknowledgements} 
The author would like to express his gratitude to Bernhard Keller for an inspiring discussion that initiated this project, as well as the sequel \cite{BorMot} (an outlook of which is given in \Cref{subsec:future}). Feedback from Xiaofa Chen, Johanne Haugland, and Calvin Pfeifer has been very useful.
Valuable input from Monica Garcia, Eric Hanson, Maximilian Kaipel, Yann Palu, and Yu Zhou is also gratefully acknowledged.

\section{Exact dg categories}\label{sec:exactdg}

\subsection{Definition} 
Recall that a $\k$-category $\c$ is \textit{additive} if 
\begin{enumerate}
\item it contains a zero object $0$, 
\item for $x,y\in\c$, one can form both the coproduct $x\sqcup y$ and product $x \sqcap y$. Furthermore, the natural morphism
$
x\sqcup y \to x \sqcap y
$
is required to be an isomorphism.
\end{enumerate}
{Since products and coproducts coincide up to isomorphism in an additive $\k$-category $\c$, we can say that $\c$ is equipped with a \textit{biproduct} $\oplus$, where $x\oplus y$ is both a product and a coproduct.}
A dg $\k$-category $\A$ is said to be \textit{additive} if $H_0\A$ is additive as a $\k$-category. Following Chen {\cite[Remark 4.5(c)]{Che24thesis}}, we will make the additional technical assumption that the $\k$-category $Z_0\A$ is additive.

\subsection{}\label{def:exact_seq}  {\cite[§3.2]{Che24thesis}} 
In order to define exact dg $\k$-categories, we will need a notion of exact sequences. 
Let $\k T$ be the path dg $\k$-category of the following dg quiver $T$ with relations:
\begin{equation*}
\begin{tikzcd}
	1 \arrow[r,"a"] & 2 \arrow[r,"b"] & 3,
\end{tikzcd}
\end{equation*}
where $a$ and $b$ are of degree 0 and have vanishing differentials, and $ba=0$. Consider the $A_{\infty}$-$k$-category $\Fun_{{A_\infty}}(\k T,\A)$ of strictly unital $A_{\infty}$-functors from $\k T$ to $\A$. We can identify $\Fun_{{A_\infty}}(\k T,\A)$ with the dg $\k$-category consisting of diagrams in $\A$ of the form
\begin{equation*}
\begin{tikzcd}
a \arrow[r,"\iota"]\arrow[rr,bend right,"\varepsilon"'] & b \arrow[r,"\pi"] & c,
\end{tikzcd}
\end{equation*}
where $\iota$ and $\pi$ are morphisms of degree 0 with vanishing differential, and $\varepsilon$ is homogeneous of degree 1 with differential $\pi\iota$. We define the \textit{category of three-term {complexes} in $\A$} by $\mathcal{H}_{3\mathrm{t}}(\A)\defeq H_0\Fun_{{A_\infty}}(\k T,\A)$. {Chen constructs a fully faithful functor $\mathcal{H}_{3\mathrm{t}}(\A) \to H_0 \Fun_{{A_\infty}}(\k\mathrm{Sq}, \A)$, where $\k\mathrm{Sq}$ is the path dg $\k$-category of the commutative square $\mathrm{Sq}$ \cite[Lemma 3.17 (and the discussion following it) and Theorem 3.16]{Che24thesis}, namely the category given by}
\begin{equation*}
\begin{tikzcd}
00\arrow[r,"\iota"]\arrow[d,"\iota'"]& 01\arrow[d,"\pi"] \\
10 \arrow[r,"\pi'"] & 11
\end{tikzcd}
\end{equation*}
{modulo the commutativity relation $\pi \iota \sim \pi' \iota'$.
A three-term complex in $\A$ is said to be a \textit{homotopy short exact sequence} if this fully faithful functor sends it to a homotopy bicartesian square \cite[Definitions 3.18 and 3.1]{Che24thesis}.}

\subsection{Definition}\label{def:exact_dg} {\cite[Definition 4.1]{Che24thesis}} 
Let $\A$ be an additive dg $\k$-category.
An \textit{exact structure} on $\A$ is given by a full subcategory $\exactstr $ of $\mathcal{H}_{3\mathrm{t}}(\A)$ {spanned by homotopy short exact sequences}, which is assumed to be closed under isomorphisms, and to obey the axioms \hyperref[item:exact_dg0]{\textbf{(Ex0)}}--\hyperref[item:exact_dg2op]{\textbf{(Ex2$\op$)}} below. An object in $\exactstr$, which is of the following form
\begin{equation}\label{eq:exact_dg_confl}
\begin{tikzcd}
a \arrow[r,tail,"\iota"]\arrow[rr,bend right,"\varepsilon"'] & b \arrow[r,two heads,"\pi"] & c,
\end{tikzcd}
\end{equation}
is called an $\exactstr$-\textit{conflation}, or just a \textit{conflation} if the choice of $\exactstr$ is natural or implicit. The map $\iota$ is called an \textit{inflation}, whereas $\pi$ is a \textit{deflation}. In our diagrams, we will decorate arrows displaying inflations and deflations as we have in \eqref{eq:exact_dg_confl}.
\begin{enumerate}
\item[\textbf{(Ex0)}]\label{item:exact_dg0} 
The identity map on a zero object in $\A$ is a deflation.\footnote{
Using \hyperref[item:exact_dg2]{\textbf{(Ex2)}} below, this axiom can be strengthened to postulate that the identity morphism on any object in $\A$ is a deflation. Since the zero object in $\mathcal{H}_{3\mathrm{t}}(\A)$ is a homotopy short exact sequence, we have that the identity map on the zero object in $\A$ is an inflation. One now uses \hyperref[item:exact_dg2op]{\textbf{(Ex2$\op$)}} to show that every identity morphism is an inflation.}
\item[\textbf{(Ex1)}]\label{item:exact_dg1} Composites of deflations are deflations.
\item[\textbf{(Ex2)}]\label{item:exact_dg2}  Deflations are stable under homotopy pullbacks in the following sense: if $\begin{tikzcd} b\arrow[r,two heads,"\pi"]&c\end{tikzcd}$ is a deflation and $\begin{tikzcd} c'\arrow[r,"f"]&c\end{tikzcd}$ is a morphism in $Z_0\A$, then the co-span
\begin{equation*}
\begin{tikzcd}
& c'\arrow[d,"f"] \\
 b \arrow[r,two heads,"\pi"] & c
\end{tikzcd}
\end{equation*}
can be completed to a homotopy pullback  
\begin{equation}\label{eq:homocarsq}
\begin{tikzcd}
b'\arrow[r,"{\pi'}",two heads]\arrow[d,"f'"]\arrow[rd,dotted]& c'\arrow[d,"f"] \\
 b \arrow[r,two heads,"\pi"] & c
\end{tikzcd}
\end{equation}
in which $\pi'$ is a deflation. For a precise definition of homotopy bicartesian squares {in $\A$}, we refer to Chen {\cite[Definition 3.25]{Che24thesis}}.
The dotted arrow in \eqref{eq:homocarsq} displays a homogeneous map of degree 1 whose differential is $\pi f' - f \pi'$.
\item[\textbf{(Ex2$\op$)}]\label{item:exact_dg2op} Inflations are stable under homotopy pushouts, in the dual sense of \hyperref[item:exact_dg2]{\textbf{(Ex2)}}.
\end{enumerate}
If $\exactstr$ is an exact structure on $\A$, we will say that the pair $(\A,\exactstr)$ (or just $\A$, if the choice of exact structure is clear) is an \textit{exact dg $\k$-category}. 

\subsection{}\label{subsec:conn} A dg $\k$-category $\A$ is called \textit{connective} if the complex $\A(a,b)$ has vanishing homology in negative\footnote{Recall from \Cref{notcon} that we are employing the homological grading convention. In the cohomological convention, connective dg $\k$-categories and $\k$-algebras are often referred to as \textit{non-positive}, since cohomology then vanishes in positive degrees.}  degrees for all $a,b\in\A$. 
A dg $k$-algebra $A$ is \textit{connective} if it is connective when regarded as a dg $\k$-category with one object.

Let $\dgcat_{\k,\geq 0}$ denote the subcategory of $\dgcat_{\k}$ spanned by the connective dg $\k$-categories. The inclusion of $\dgcat_{\k,\geq 0}$ into $\dgcat_{\k}$ admits a right adjoint $\tau_{\geq 0}$. For $\A\in \dgcat_{\k}$, consider the co-unit $\tau_{\geq 0}\A \to \A$ of the adjunction. This is the \textit{connective cover} of $\A$. The connective cover induces a bijection of exact structures of exact structures on $\A$ and exact structures on $\tau_{\geq 0}\A$ \cite[Remark 4.25]{Che24thesis}. Note also that the induced functor $H_0(\tau_{\geq 0}\A) \to H_0(\A) $ is an equality of $\k$-categories.

\subsection{Definition} 
Let $(\A,\exactstr)$ and $(\A',\exactstr')$ be exact dg $\k$-categories, and let $F\colon \A\to \A'$ be a morphism in $\Hqe_{\k}$. If the induced functor $\mathcal{H}_{3\mathrm{t}}(F)\colon \mathcal{H}_{3\mathrm{t}}(\A)\to \mathcal{H}_{3\mathrm{t}}(\A')$ restricts to a functor $\mathcal{H}_{3\mathrm{t}}(F)|_{\exactstr}\colon \exactstr\to \exactstr'$, we will call $F$ an \textit{exact morphism} in $\Hqe_{\k}$. 
Let $\Hqe_{\k}^{\mathrm{ex}}({(\A,\exactstr),(\A',\exactstr')})$ denote the subset of $\Hqe_{\k}(\A,\A')$ consisting of exact morphisms. {We abbreviate the notation to $\Hqe_{\k}^{\mathrm{ex}}(\A,\A')$ whenever the exact structures on $\A$ and $\A'$ are unequivocally imposed by the context.}

\subsection{}
Let $\mathscr{C}(\k)$ denote the (typically large) dg $\k$-category of complexes of $\k$-modules, and for an additive dg $\k$-category $\A$, let $\mathscr{C}(\A)$ denote the (also typically large) dg $\k$-category of dg functors $\A\op\to \mathscr{C}(\k)$. There is a natural way to equip $\mathscr{C}(\k)$ with an exact structure, which in turn induces an exact structure on $\mathscr{C}(\A)$. Moreover, we can define a \textit{shift} of an object $X\in \mathscr{C}(\A)$ by degree-shifting complexes in $\mathscr{C}(\k)$. 
The dg $\k$-category $\A$ is said to be \textit{pretriangulated} if the image of $\A$ under the Yoneda embedding
\begin{equation*}
	\begin{tikzcd}
		\A \arrow[r,hook] & \mathscr{C}(\A)
	\end{tikzcd}
\end{equation*}
is closed under extensions and shifts \cite{BK90}. Since morphisms in a pretriangulated dg $\k$-category admit a homotopy kernel and a homotopy cokernel, and the class of conflations that arise in this way define an exact structure on $\A$ \cite[Example 4.7]{Che24thesis}. Pretriangulated dg $\k$-categories are thus examples of exact dg $\k$-categories. Unless otherwise specified, we will denote the suspension endofunctor on a pretriangulated dg $\k$-category by $\Sigma$. This will also be the case for triangulated $\k$-categories.

\subsection{} \label{Che23.5.7} {\cite[Example-Definition 4.8]{Che24thesis}} 
Let $\mathfrak{a}$ be a full dg subcategory of an exact dg $\k$-category $\A$. We say that $\mathfrak{a}$ is \textit{extension-closed in $\A$} provided that for all conflations in $\A$
\begin{equation*}
\begin{tikzcd}
a \arrow[r,tail,""]\arrow[rr,bend right,""'] & b \arrow[r,two heads,""] & c
\end{tikzcd}
\end{equation*}
with $a,c\in\mathfrak{a}$, we also have $b\in \mathfrak{a}$. In this case, the conflations in $\A$ with terms in $\mathfrak{a}$ {yield} an exact structure on $\mathfrak{a}$ in such a way that the inclusion functor $\mathfrak{a}\embed \A$ is an exact morphism in $\Hqe_{\k}$. {This is the maximal exact structure with which $\mathfrak{a}$ can be equipped so that this morphism becomes exact.}

\subsection{} \label{def:Db}
Let $\A$ be a small connective exact dg $\k$-category.
The \textit{bounded derived dg $\k$-category} of $\A$ is denoted $\Ddgb(\A)$, and is defined as follows: the exact morphism $\A\to \Ddgb(\A)$ in $\Hqe_{\k}^{\mathrm{ex}}$ is to satisfy the universal property that for every pretriangulated dg $\k$-category $\B$, the induced map
\begin{equation}\label{eq:UVP_Ddgb}
	\Hqe_{\k}^{\mathrm{ex}}(\Ddgb(\A),\B) \to \Hqe_{\k}^{\mathrm{ex}}(\A,\B) 
\end{equation}
is a bijection. This universal morphism indeed exists \cite[Theorem 6.1]{Che23thesis}. 
\subsection{} 
The bounded derived dg $\k$-category of a small connective exact dg $\k$-category is small by construction. For an arbitrary additive dg $\k$-category $\mathfrak{a}$, we will also need to consider the \textit{derived $\k$-category} $\Der(\mathfrak{a})$ of $\mathfrak{a}$, which is defined as the localization of $\mathscr{C}(\mathfrak{a})$ with respect to the \textit{quasi-isomorphisms}, namely the class of morphisms in $\mathscr{C}(\mathfrak{a})$ inducing quasi-isomorphisms of complexes when evaluated at any object in $\mathfrak{a}$. The derived $\k$-category of $\mathfrak{a}$ is pretriangulated, but not hardly ever small.

\subsection{Example}\label{eg:exactdg}
\begin{enumerate}
\item {\cite[Example 4.6]{Che24thesis}} 
Any $\k$-category can be regarded as a dg $\k$-category in which all mapping complexes are stalk complexes. 
Specializing the axioms of exact dg $\k$-categories to $\k$-categories recovers the notion of Quillen exact $\k$-categories \cite{Qui73}  \cite[Appendix A]{Kel90}. Consequently, any Quillen exact $\k$-category can be regarded as an exact dg $\k$-category. In particular, any abelian $\k$-category is naturally equipped with an exact structure, if regarded as a dg $\k$-category.
\item\label{eg:exactdg1} Inside the derived dg $\k$-category $\Der(\mathfrak{a})$ of an additive dg $\k$-category $\mathfrak{a}$ sits the smallest pretriangulated subcategory {closed under direct summands} containing $\mathfrak{a}$. This is denoted $\per(\mathfrak{a})$, and called the \textit{perfect derived dg $\k$-category} of $\mathfrak{a}$. We note that although $\Der(\mathfrak{a})$ is not small, the dg subcategory $\per(\mathfrak{a})$ is. Inside $\per(\mathfrak{a})$ sits the \textit{two-term dg $\k$-category of $\mathfrak{a}$}, namely the full subcategory spanned by the objects $x$ that appear in a conflation of the form
\begin{equation*}
\begin{tikzcd}
a_0 \arrow[r,tail,""]\arrow[rr,bend right,""'] & x \arrow[r,two heads,""] & \Sigma a_1
\end{tikzcd}
\end{equation*}
where $a_0$ and $a_1$ are in the essential image of the embedding $\begin{tikzcd} \mathfrak{a} \arrow[r,hook] & \per(\mathfrak{a})\end{tikzcd}\!\!$. One sees that this defines an extension-closed dg subcategory of $\per(\mathfrak{a})$, whence it admits the structure of exact dg $\k$-category by \Cref{Che23.5.7}.
Throughout this text, the two-term dg $\k$-category of $\mathfrak{a}$ will be denoted $\per_{[0,1]}(\mathfrak{a})$. When $\mathfrak{a}$ consists of a single object whose endomorphism dg $\k$-algebra is $A$, we write $\per_{[0,1]}(A)$ instead of $\per_{[0,1]}(\mathfrak{a})$. The \textit{two-term category} of $\mathfrak{a}$ is defined as $H_0\per_{[0,1]}(\mathfrak{a})$.
\end{enumerate}

\section{Extriangulated categories}\label{subsec:algextri}

 \subsection{}
In order to recall the definition of an extriangulated $\k$-category, some preparation of terminology will be needed. Although the original axioms are due to Nakaoka--Palu \cite[Definition~2.12]{NP19}, we will present the axioms of 1-exangulated $\k$-categories, an equivalent notion put forward by Herschend--Liu--Nakaoka \cite[Definition 2.32 and Proposition 4.3]{HLN17}.

\subsection{}\label{def:Eext} Fix an additive $\k$-category $\c$. Given a $\k$-linear bifunctor 
\begin{equation}\label{eq:Eext_bifun}
    \mathbb{E}\colon \c\op\times\c \to \Mod(\k),
\end{equation}
we will refer to an element $\xi\in \E(c,a)$ as \textit{$\E$-extension}, where the objects $a,c\in \c$ are arbitrary. A \textit{morphism of $\E$-extensions} from $\xi\in\E(c,a)$ to $\xi'\in\E(c',a')$ is a given pair $(f\colon a \to a',g\colon c\to c')$ of morphisms in $\c$ such that $\E(g,a')(\xi')=\E(c,f)(\xi)$. Let $\E\mhyphen \mathrm{Ext}$ denote the $\k$-category of $\E$-extensions.

We can regard $\c$ as a dg $\k$-category with trivial differential grading, and then define $\mathcal{H}_{3\mathrm{t}}(\c)$ as in \Cref{def:exact_seq}, which becomes a category of three-term {complexes} in the conventional sense. Explicitly, the category $\mathcal{H}_{3\mathrm{t}}(\c)$ is that of representations in $\c$ of the quiver 
\begin{equation*}
\begin{tikzcd}
1 \arrow[r] & 2 \arrow[r] & 3,
\end{tikzcd}
\end{equation*}
bound by the relation annulling the path from $1$ to $3$.
For a $\k$-linear bifunctor $\E$ as displayed in \eqref{eq:Eext_bifun}, we define an \textit{exact realization} of $\E$ as a functor
\begin{equation}\label{eq:realiz_fun}
	\mathfrak{s}\colon \E\mhyphen \mathrm{Ext} \to \mathcal{H}^{}_{3\mathrm{t}}(\c)
\end{equation}
with the following two properties:
\begin{enumerate}
	\item\label{realiz1} the $\E$-extension $0\in \E(c,a)$ is sent to a \textit{split sequence} $a \to a\oplus c\to c$, i.e. the first morphism is natural when $a\oplus c$ is regarded as a coproduct of $a$ and $c$, and the latter is natural when $a\oplus c$ is regarded as a product.
	\item For  $\xi\in \E(c,a)$ and $\xi'\in \E(c',a')$, we have that $\mathfrak{s}(\xi\oplus \xi') = \mathfrak{s}(\xi)\oplus \mathfrak{s}(\xi')$ in $\E(c\oplus c',a\oplus a')$.
	 \end{enumerate}
A definition along these lines has been put forward in previous work of the author with P. Trygsland \cite[Definition 4.4]{BT21}.
Consider the subcategory $\mathcal{H}^{y}_{3\mathrm{t}}(\c)$ of $\mathcal{H}_{3\mathrm{t}}(\c)$ containing all objects, but containing only the morphisms of the form 
\begin{equation*}
\begin{tikzcd}
a \arrow[r,"\iota"]\arrow[d,equal] & b\arrow[r,"\pi"]\arrow[d] & c\arrow[d,equal] \\ 
a \arrow[r,"\iota'"]  & b' \arrow[r,"\pi'"] & c
\end{tikzcd}
\end{equation*}
where the outer morphisms are identity morphisms in $\c$. 
For each pair of objects $a,c\in\c$, the functor $\mathfrak{s}$ induces a map from $\E(c,a)$ to the set of isomorphism classes in $\mathcal{H}^{y}_{3\mathrm{t}}(\c)$ of conflations of the form
\begin{equation}\label{eq:realiz_conf}
\begin{tikzcd}
a \arrow[r,"\iota",tail]& b\arrow[r,"\pi",two heads]& c.
\end{tikzcd}
\end{equation}
Our functorial definition of exact realizations is thus consistent with the definitions elsewhere in the literature \cite[§2.2]{NP19} \cite[§2.4]{HLN17}. 

Given an {exact realization} $\mathfrak{s}$ of $\E$, consider the three-term {complex} in \eqref{eq:realiz_conf} as an object in $\mathcal{H}_{3\mathrm{t}}(\c)$.
The three-term {complexes} in the image of $\mathfrak{s}$ will be called \textit{$\mathfrak{s}$-conflations}, or just \textit{conflations} when $\mathfrak{s}$ can be safely suppressed from the notation. A \textit{morphism of ($\mathfrak{s}$-)conflations} is simply a morphism of the underlying complexes. The morphism $\iota$ in \eqref{eq:realiz_conf} is called an \textit{($\mathfrak{s}$-)inflation}, whereas $\pi$ is a(n \textit{$\mathfrak{s}$-)deflation}. In our diagrams, we will decorate inflations and deflations as we have in \eqref{eq:realiz_conf}.

\subsection{Definition}\label{def:extri} \cite[Definition 2.32]{HLN17}. Let $\c$ be an additive $\k$-category, let
\begin{equation*}
    \mathbb{E}\colon \c\op\times\c \to \Mod(\k),
\end{equation*}
be an additive bifunctor, and let $\mathfrak{s}$ be an exact realization of $\E$. Then the triple $(\c,\E,\mathfrak{s})$ is called an \textit{extriangulated $\k$-category} if the following axioms hold:
\begin{enumerate}
\item[\textbf{(Extr1)}]\label{item:extri1} Composites of $\mathfrak{s}$-deflations are $\mathfrak{s}$-deflations, and composites of $\mathfrak{s}$-inflations are $\mathfrak{s}$-inflations.
\item[\textbf{(Extr2)}]\label{item:extri2} $\mathfrak{s}$-deflations admit \textit{good lifts} in the following precise sense: given the solid part of the diagram \begin{equation*}
\begin{tikzcd}
a \arrow[r,tail,"\iota'",dashed]\arrow[d,dashed,"1_a"] & b'\arrow[r,two heads,"\pi'",dashed]\arrow[d,dashed,"f'"] & c'\arrow[d,"f"] \\
 a \arrow[r,tail,"\iota"]  & b \arrow[r,two heads,"\pi"] & c
\end{tikzcd}
\end{equation*}
where the bottom row is an $\mathfrak{s}$-conflation, we can add the dashed arrows to produce a morphism of $\mathfrak{s}$-conflations in which the leftmost vertical morphism is an identity morphism. Moreover, we require the ``folded'' complex 
\begin{equation*}
\begin{tikzcd}[ampersand replacement=\&, column sep=4em]
b' \arrow[r,"{\begin{pmatrix} f' \\ \pi' \end{pmatrix}}"] \& b\oplus c' \arrow[r,"{\begin{pmatrix} -\pi &  f \end{pmatrix}}"] \& c
\end{tikzcd}
\end{equation*}
to be an $\mathfrak{s}$-conflation.
\item[\textbf{(Extr2$\op$)}]\label{item:extri2op} $\mathfrak{s}$-inflations admit good lifts, in the dual sense of \hyperref[item:extri2]{\textbf{(Extr2)}}.
\end{enumerate}
If $(\c,\E,\mathfrak{s})$ is an {extriangulated $\k$-category}, we say that that $\E$ and $\mathfrak{s}$ form an \textit{extriangulation} on the underlying additive $\k$-category $\c$.
The extriangulated $\k$-category $(\c,\E,\mathfrak{s})$ will be denoted $(\c,\E)$, or simply $\c$, whenever reference to the exact realization and/or the bifunctor is unnecessary.

\subsection{} Given a conflation
\begin{equation*}
\begin{tikzcd}
a \arrow[r,"\iota",tail]& b\arrow[r,"\pi",two heads]& c
\end{tikzcd}
\end{equation*}
in an extriangulated $\k$-category $\c$, we say that
\begin{enumerate}
	 \item $c$ is a \textit{cone} of $\iota$,
	 \item $a$ is a \textit{cocone} of $\pi$,
	 \item $b$ is an \textit{extension} of $c$ in $a$.
\end{enumerate}
Analogously to \Cref{Che23.5.7}, we say that a full subcategory $\x\subseteq \c$ is \textit{extension-closed} if it is closed under the formation of extensions. We define the notions of \textit{closure under cones} and \textit{closure under cocones} similiarly.

\subsection{Definition}\label{def:extrifun} \cite[Definition 2.32]{B-TS20} \cite[Definition 4.6]{BT21} 
Let $(\c,\E,\mathfrak{s})$ and $(\c',\E',\mathfrak{s}')$ be extriangulated $\k$-categories. A $\k$-linear \textit{extriangulated functor} from $(\c,\E,\mathfrak{s})$ to $(\c',\E',\mathfrak{s}')$ is given by a pair $(F,\psi)$ where $F\colon \c\to \c'$ is a $\k$-linear functor and $\psi\colon \E(-,-) \Longrightarrow \E'(F-,F-)$ is a natural transformation of functors $\c\op\times\c \to \Mod(\k)$ such that the square in \eqref{eq:extrifun} below commutes up to natural isomorphism. In this square, the functor $F_{\psi} \colon \E\mhyphen \mathrm{Ext} \to  \E'\mhyphen \mathrm{Ext}$ is induced by $F$ and $\psi$ in the sense that $\xi\in \E\mhyphen \mathrm{Ext}$ is sent to $\psi \xi \in \E'\mhyphen \mathrm{Ext}$. 
\begin{equation}\label{eq:extrifun}
	\begin{tikzcd}[column sep=4em]
		 \E\mhyphen \mathrm{Ext} \arrow[r,"F_{\psi}"]\arrow[d,"\mathfrak{s}"]  &  \E'\mhyphen \mathrm{Ext}\arrow[d,"\mathfrak{s}'"] \\ 
		 \mathcal{H}^{}_{3\mathrm{t}}(\c) \arrow[r,"\mathcal{H}^{}_{3\mathrm{t}}(F)"]  &  \mathcal{H}^{}_{3\mathrm{t}}(\c')
	\end{tikzcd}
\end{equation}
Whenever it is unnecessary to explicitly mention the natural transformation giving $F$ the structure of an extriangulated functor, we will refrain from doing so. We will often say that a $\k$-linear functor $F$ \textit{is extriangulated} when there is an obvious choice of natural transformation $\psi$ such that the pair $(F,\psi)$ is an extriangulated functor.

A $\k$-linear extriangulated functor $(F,\eta)\colon \c \to \c'$ is a \textit{$\k$-linear extriangulated equivalence} if the functor $F$ is an equivalence of $\k$-categories and $\eta$ is a natural isomorphism \cite[Proposition 2.13]{NOS22}. It is equivalent that there exists an extriangulated functor $(G,\varepsilon)\colon \c' \to \c$ such that $G$ is a quasi-inverse of $F$ and the composite natural isomorphisms $\varepsilon \eta$ and $\eta \varepsilon$ are the identity transformations on $\E$ and $\E'$, respectively \cite[Proposition 4.11]{BTHSS23}.

An extriangulated $\k$-category $(\c',\E',\mathfrak{s}')$ is an \textit{extriangulated subcategory} of an extriangulated $\k$-category $(\c,\E,\mathfrak{s})$ if $\c'$ is a subcategory $\c$, the inclusion functor $I\colon\c' \embed \c$ induces a natural monomorphism $\eta$ from $\E'$ to the restriction of $\E$ to the essential image of $\c$, and the pair $(I,\eta)$ is an extriangulated functor. 

\subsection{}\label{lem:NP19.2.18}\cite[Remark 2.18]{NP19}
We have the following counterpart of \Cref{Che23.5.7}:
Let $\c'$ be an extension-closed subcategory of an extriangulated $\k$-category $(\c,\E)$. 
The restriction $\E'$ of $\E$ to $\c$ determines an extriangulation on $\c'$, and $(\c',\E')$ is an extriangulated subcategory of $(\c,\E)$.

\subsection{}\label{rem:extri_eg}
The structure of a triangulated $\k$-category is a special instance of extriangulated structure \cite[Example 2.13]{NP19}, and triangulated functors are examples of extriangulated functors \cite[Theorem 2.33]{B-TS20}. 
More generally, if $(\c,\E,\mathfrak{s})$ is an extension-closed subcategory of a triangulated $\k$-category $(\mathcal{D},\mathbb{F},\mathfrak{t})$ and we have a triangulated functor $(\mathcal{D},\mathbb{F},\mathfrak{t})\to^{F} (\mathcal{D}',\mathbb{F}',\mathfrak{t}')$, then $F$ restricts to an extriangulated functor $(\c,\E,\mathfrak{s})\to^{F|_{\c}} (\mathcal{D}',\mathbb{F}',\mathfrak{t}')$. Should we replace the codomain of $F|_{\c}$ with any extriangulated subcategory of $(\mathcal{D}',\mathbb{F}',\mathfrak{t}')$ containing the essential image of the composite 
$$\begin{tikzcd} \c \arrow[r,hook] & \mathcal{D}  \arrow[r,"F"] & \mathcal{D}',\end{tikzcd}$$ we still obtain an extriangulated functor.

\subsection{}\label{Che23.6.19}
Let $\A$ be an exact dg $\k$-category.
Then the $\k$-category $H_0\A$ can be equipped with an extriangulation \cite[Theorem 4.26]{Che24thesis}. This generalizes Bondal--Kapranov's result on the triangulation on the zeroth homology category of pretriangulated dg $\k$-category \cite[Proposition 2]{BK90}. Whenever the exact structure $\exactstr$ is suppressed from the notation, we will denote the bifunctor in \eqref{eq:Che23.6.19_E} below by $\E_{\A}$.  In the same spirit, we denote the associated extriangulated $\k$-category by the pair by $(H_0\A,\E_{\A})$, or simply by $H_0\A$. An extriangulated $\k$-category is called \textit{algebraic} if there exists a $\k$-linear extriangulated equivalence to an extriangulated $\k$-category of the form $(H_0\A,\E_{\A})$ for some exact dg $\k$-category $\A$.  
If an extriangulated $\k$-category $(\c,\E)$ is equivalent to $(H_0\A,\E_{\A})$, we will say that $(\c,\E)$ is \textit{enhanced} by $\A$, or that $\A$ is a(n \textit{exact) dg enhancement} of $(\c,\E)$.

We briefly recall how the extriangulation on $H_0\A$ comes about. Chen constructs the bifunctor 
\begin{equation}\label{eq:Che23.6.19_E}
 \E_{\A}\colon (H_0\A)\op\times H_0\A \to \Mod(\k),
 \end{equation}
in such a way that $\E_{\A}(c,a)$ is the set of isomorphism classes in $\mathcal{H}^{y}_{3\mathrm{t}}(\c)$ of objects of the form \eqref{eq:realiz_conf} {\cite[§4.2]{Che24thesis}}. One can show that the category of $\E_{\A}$-extensions becomes isomorphic to $\mathcal{H}^{}_{3\mathrm{t}}(\c)$, and that this equivalence of categories serves as a realization of $\E_{\A}$. 

\subsection{}\label{enhtri}
Suppose that $\A$ is an exact dg $\k$-category, and let $\tau_{\geq 0}\A\to \A$ be its connective cover (as defined in \Cref{subsec:conn}). 
The universal morphism $\tau_{\geq 0}\A\to \Ddgb(\tau_{\geq 0}\A)$ in $\Hqe_{\k}$ (see \Cref{def:Db}) induces a fully faithful extriangulated functor 
\begin{equation}\label{eq:H0enhtri}
\begin{tikzcd}
H_0 \A  = H_0(\tau_{\geq 0}\A)\arrow[r,hook]  & H_0 \Ddgb(\tau_{\geq 0}\A).
\end{tikzcd} 
\end{equation}
Since $\Ddgb(\tau_{\geq 0}\A)$ is a pretriangulated dg $\k$-category, we have that $H_0\Ddgb(\tau_{\geq 0}\A)$ is an algebraic triangulated $\k$-category \cite[2.1, 2.2]{Kel99}. It follows that any extriangulated $\k$-category admitting a dg enhancement can be obtained as an extension-closed subcategory of an algebraic triangulated $\k$-category.

\subsection{} Let $\A$ and $\A'$ be exact dg $\k$-categories and let $F\colon \A \to \A'$ be an exact morphism in $\Hqe_{\k}$. By the universal property of the bounded derived dg $\k$-category, we can induce an exact morphism $\widetilde{F}\colon \Ddgb(\A) \to \Ddgb(\A')$ in such a way that we obtain the following commutative square in $\Hqe_{\k}$:
\begin{equation*}
\begin{tikzcd}
	\A \arrow[r,"F"]\arrow[d] & \A'\arrow[d] \\ 
	\Ddgb(\tau_{\geq 0}\A) \arrow[r,"\widetilde{F}"] & \Ddgb(\tau_{\geq 0}\A')
\end{tikzcd}
\end{equation*}
It follows from \Cref{rem:extri_eg} that the functor $H_0F\colon H_0\A \to H_0\A'$ can be equipped with a natural transformation so that it becomes an extriangulated functor. 

\section{Localization}\label{subsec:loc}

\subsection{} \cite{GZ12} Let $\c$ be a category and let $\mathcal{S}$ be a set of morphisms in $\c$. A \textit{Gabriel--Zisman localization} of $\c$ with respect to $\mathcal{S}$ is a functor $L_\mathcal{S} \colon \c \to \c[\mathcal{S}^{-1}]$ having the following properties:
\begin{enumerate}
	\item Any morphism $s\in \mathcal{S}$ is sent to an isomorphism by $L_{\mathcal{S}}$.
	\item Let $F\colon \c\to \c'$ be a functor sending all morphisms in $\mathcal{S}$ to isomorphisms. Then, up to natural isomorphism of functors, the functor $F$ factors uniquely through $L_{\mathcal{S}}$. 
\end{enumerate}
A Gabriel--Zisman localization always exists (up to set-theoretic issues), and is clearly unique up to natural isomorphism. In the cases where $\c$ is closed under composition, one may also refer to $\mathcal{S}$ as a subcategory of $\c$.

\subsection{}\label{idquot} Let $\c$ be a $\k$-category and let $\calr$ be a full subcategory of $\c$ with is closed under direct summands. 
Letting $[\calr]$ is the two-sided ideal of $\c$ consisting of morphisms factoring through an object in $\calr$, one defines the \textit{ideal quotient} of $\c$ with respect to $\calr$ by $\c/[\calr]$. The natural functor $\overline{(-)}\colon \c \to \c/[\calr]$ will be called the \textit{ideal quotient functor}. One can attain the ideal quotient functor as a Gabriel--Zisman localization by considering the class of morphisms in $\c$ consisting of split monomorphisms admitting a cokernel in $\calr$ \cite[Example 2.6]{Oga22}. If $(\c,\E)$ is an extriangulated $\k$-category such that $\E(x,r) = 0 =\E(r,x)$ for all $r\in \calr$ and $x\in\c$, one naturally induces a $\k$-linear bifunctor $\overline{\E}$ on $\c/[\calr]$, and thus an extriangulation on $\c/[\calr]$ \cite[Proposition 3.30]{NP19}. The ideal quotient functor $\overline{(-)}\colon \c \to \c/[\calr]$ becomes an extriangulated functor when equipped with the natural epimorphism $\E\Longrightarrow \overline{\E}$.

\subsection{}\label{NOSloc}
Let $(\c,\E)$ be an extriangulated $\k$-category. Alike the special case of triangulated $\k$-categories \cite[§2.1]{Ver96}, one can use Gabriel--Zisman localization of develop a notion of localization.
A \textit{thick subcategory} of $(\c,\E)$ is a full subcategory $\mathcal{N}\subseteq \c$ which is
\begin{enumerate}
	\item closed under isomorphisms and direct summands,
	\item it is extension-closed, closed under cones, and closed under cocones.
\end{enumerate}
The smallest thick subcategory of $\c$ containing a full subcategory $\x\subseteq \c$ will be denoted $\thick_{\c}(\x)$, or simply $\thick(\x)$ when the ambient extriangulated $\k$-category $\c$ is clear from context. For a thick subcategory $\mathcal{N}$ of $(\c,\E)$ let $\mathit{Inf}_{\mathcal{N}}$ be the set of inflations in $\c$ with cone in $\mathcal{N}$. Dually, we define $\mathit{Def}_{\mathcal{N}}$ as the set of deflations in $\c$ with cocone in $\mathcal{N}$. Defining the subcategory $\mathcal{S}_{\mathcal{N}}$ of $\c$ as the smallest subcategory of $\c$ containing both $\mathit{Inf}_{\mathcal{N}}$ or $\mathit{Def}_{\mathcal{N}}$, one defines the \textit{Verdier localization functor} of $\c$ with respect to $\mathcal{N}$ as the Gabriel--Zisman localization of $\c$ with respect to $\mathcal{S}_{\mathcal{N}}$. 
It will be denoted $L_{\mathcal{N}}\colon(\c,\E)\to (\c/\mathcal{N},{\E_{\mathcal{N}}})$. The codomain $\c/\mathcal{N}$ of this functor is called the \textit{Verdier quotient} of $\c$ with respect to $\mathcal{N}$. We say that the thick subcategory $\mathcal{N}$ \textit{admits a calculus of left fractions} if the subcategory $\mathcal{S}_{\mathcal{N}}$ admits a calculus of left fractions, in the sense of Gabriel--Zisman \cite[§I.2]{GZ12}.
It is not true in general that the Verdier quotient $\c/\mathcal{N}$ admits an extriangulation, even when $\mathcal{S}_{\mathcal{N}}$ admits a calculus of left fractions; the axiom \hyperref[item:extri1]{\textbf{(Extr1)}} need not hold unless $\mathcal{N}$ is biresolving or percolating \cite[§4]{NOS22}.
However, since $\mathcal{S}_{\mathcal{N}}$ contains all isomorphisms and is closed under the composition and the formation of direct sums, the Verdier quotient $\c/\mathcal{N}$ becomes weakly extriangulated, in the sense of Baillargeon--Br{\"u}stle--Gorsky--Hassoun \cite{BBGH22}.

\subsection{}\label{dgloc}
Let $\A$ be a dg $\k$-category and let $\mathcal{S}$ be a subcategory of $H_0\A$. For a dg $\k$-category $\B$, let $\Hqe_{k}(\A,\B)_{\mathcal{S}}$ denote the subset of $\Hqe_{\k}(\A,\B)$ consisting of the exact morphisms $\A\to \B$ inducing functors $H_0\A\to H_0\B$ that send every $s\in\mathcal{S}$ to an isomorphism.
The \textit{localization morphism} $\varpi_{\mathcal{S}}\colon\A\to \L_{\mathcal{S}}\A$ is defined in terms of a universal property; this is the morphism in $\Hqe_{\k}$ such that for every $\B\in\dgcat_{\k}$, the induced map 
\begin{equation}\label{eq:UVP_loc}
	\varpi_{\mathcal{S}}^{\ast}\colon \Hqe_{k}(\L_{\mathcal{S}}\A,\B) \to \Hqe_{k}(\A,\B)
\end{equation}
is an injection whose image is $\Hqe_{k}(\A,\B)_{\mathcal{S}}$ \cite{Toe07}. The dg $\k$-category $\L_{\mathcal{S}}\A$, or more precisely any dg $\k$-category to which it is isomorphic in $\Hqe_{\k}$, is called a \textit{dg localization} of $\A$ with respect to $\mathcal{S}$. 

\subsection{} 
A dg $\k$-category $\A$ \textit{has $\k$-cofibrant mapping complexes} if any mapping complex $\A(x,y)$ is $\k$-cofibrant in the projective model structure on the category of chain complexes of $\k$-modules.
If $\A$ has $\k$-cofibrant mapping complexes, its dg localization takes a concrete form. For a full rigid dg subcategory $\r\subseteq \A$, we define the \textit{dg quotient} $\A/\r$ by freely adding a morphism $\varepsilon_p\colon p\to p$ of degree 1 with $\partial \varepsilon_p = 1_p$, for every $p\in\r$ \cite{Dri04}. The {dg quotient morphism} $\A \to \A/\r$ satisfies the following universal property \cite[Theorem 4.0.1]{Tab10}: this is the morphism in $\Hqe_{\k}$ such that for every $\B\in\dgcat_{\k}$, the induced map 
\begin{equation}\label{eq:UVP_loc_drin}
	\varpi_{\r}^{\ast}\colon \Hqe_{k}(\A/\r,\B) \to \Hqe_{k}(\A,\B)
\end{equation}
is an injection whose image is the subset of $\Hqe_{\k}(\A,\B)$ consisting of the exact morphisms $\A\to \B$ inducing functors $H_0\A\to H_0\B$ sending $H_0\r$ to zero objects.

  
\subsection{Lemma}\label{prop:dgquot_ex}
\begin{enumerate}
\item\label{prop:dgquot_ex1} Let $\A$ be a connective additive dg $\k$-category and let $\mathcal{S}$ be a subcategory of $H_0\A$. Suppose that the localized $\k$-category $H_0(\A)[\mathcal{S}^{-1}]$ is additive and that the localization functor $H_0(\A)\to H_0(\A)[\mathcal{S}^{-1}]$ is an additive functor. We then have an equivalence of $\k$-categories $H_0(L_{\mathcal{S}}\A) \simeq H_0\A[\mathcal{S}^{-1}]$.
\item\label{prop:dgquot_ex2} If $\A$ is an exact dg $\k$-category with $\k$-cofibrant mapping complexes and $\r$ is a full dg subcategory of $\A$, we have equivalences
\begin{equation*}
	 H_0(\A/\r) \simeq H_0(L_{\mathcal{S}_{\thick(\calr)}}\A) \simeq H_0(\A)/\thick(H_0\r).
\end{equation*}
\end{enumerate}
\begin{proof}
	Chen proves \ref{prop:dgquot_ex1} in the case where $\mathcal{S}$ admits a calculus of left fractions \cite[Lemma 5.37]{Che23thesis}. The proof uses the fact that the localization functor is additive,\footnote{This is the case when $\mathcal{S}$ admits a left calculus of fractions \cite[I.3.3]{GZ12} \cite[\href{https://stacks.math.columbia.edu/tag/05QE}{Tag 05QE}]{stacks}.} and then exploits universal properties. We can thus upgrade the statement as we have above, and assert that \ref{prop:dgquot_ex1} holds.\footnote{In private conversation, X. Chen convinced the author that is argument is valid. This discussion is gratefully acknowledged.} To prove \ref{prop:dgquot_ex2}, one first argues that $\A/\r$ satisfies the universal property in \eqref{eq:UVP_loc} with $\mathcal{S} = \mathcal{S}_{\thick(H_0\r)}$, as defined in \Cref{NOSloc}. Next, we first appeal to \Cref{subsec:conn} in order to freely replace $\A$ with its connective cover $\tau_{\geq 0}\A$. Since $H_0\tau_{\geq 0}\A=H_0\A$, the result now follows from \ref{prop:dgquot_ex1}.
	\end{proof}

 \section{Silting in (0-Auslander) extriangulated categories}\label{subsec:0Ausdef}
 
 \subsection{}\label{def:proj_inj}
An object $x$ in an extriangulated $\k$-category $(\c,\E)$ is \textit{$\E$-projective} (resp. \textit{$\E$-injective}) if the functor $\E(x,-)\colon \c\to \Mod(\k)$ (resp. $\E(-,x)\colon \c\op\to \Mod(\k)$) is zero. An \textit{$\E$-projective-injective object} in $\c$ is both $\E$-projective and $\E$-injective. We will simply use the terms \textit{projective}, \textit{injective} and \textit{projective-injective} when there can be no confusion as to the choice of bifunctor $\E$.
Let $\proj(\c)$ (resp. $\inj(\c)$, resp. $\projinj(\c)$) denote the full subcategory of $\c$ spanned by the projective (resp. injective, resp. projective-injective) objects in $\c$. 
An extriangulated $\k$-category $(\c,\E)$ \textit{has enough projectives} if one for every object $x\in \c$ can find a deflation $\begin{tikzcd} q \arrow[r,two heads,"\pi"] & x\end{tikzcd}\!\!$, where $q$ is projective. Such a deflation will be called a \textit{projective cover} of $x$.
The notions of \textit{having enough injectives} and \textit{injective envelopes} are defined dually. 

 \subsection{}\label{defprop:GNP23.2.1} The extriangulated $\k$-category $(\c,\E)$ is \textit{hereditary} if all functors of the form $\E(x,-)\colon \c\to \Mod(\k)$ are right exact. It is equivalent to impose the right exactness of all contravariant functors of the form $\E(-,y)$. 
In the case where $(\c,\E)$ has enough projectives, it is hereditary if and only if one for all objects $x\in\c$ can find a conflation $\begin{tikzcd} p_1\arrow[r,tail] &  p_0 \arrow[r,two heads] & x\end{tikzcd}$ in which $p_0$ and $p_1$ are projective \cite[Proposition 2.1]{GNP23}. One formulates a dual condition when $(\c,\E)$ has enough injectives. 

 \subsection{}\label{rem:her_highext} 
In the same spririt as Yoneda's construction for abelian $\k$-categories \cite{Yon60}, one defines the higher extension functors $\E^i$, for $i\geq 0$, of $(\c,\E)$ in terms of coend powers of $\E$ \cite[Definition 3.1]{GNP21}. The functor $\E^0$ is the Hom-bifunctor on $\c$, whereas $\E^1 = \E$, both by definition.
 If $(\c,\E)$ is algebraic and enhanced by $\A$, one can recover the higher extension functors $\E^i$ as higher Ext-functors in the bounded derived dg $\k$-category of $\A$ \cite[Proposition 6.24]{Che23thesis}.

{ Suppose that the extriangulated $\k$-category $(\c,\E)$ has enough projectives. For an object $x\in \c$, consider a conflation }
\begin{equation}\label{eq:LN19.5.2}
\begin{tikzcd}
x'\arrow[r,tail] & p \arrow[r,two heads] & x,
\end{tikzcd}
\end{equation}
{where $p$ is projective. Then, for any $i\geq 2$, we have that $\E^i(x,-)$ and $\E^{i-1}(x',-)$ are naturally isomorphic as covariant functors from $\c$ to $\mathrm{Mod}(\k)$ \cite[Proposition 5.2]{LN19}. A dual result holds when $(\c,\E)$ has enough injectives.}

\subsection{}\label{les} {Let }
$\begin{tikzcd}
a\arrow[r,tail] & b \arrow[r,two heads] & c
\end{tikzcd}$
{be a conflation in an extriangulated $\k$-category $(\c,\E)$ and let $x\in \c$. We then have long-exact sequences \cite[Proposition 5.2]{LN19}}
\begin{equation}\label{eq:rightles}
	\begin{tikzcd}[row sep=0.95em]
                    {} & {\c(c,x)} \arrow[r]   & {\c(b,x)} \arrow[r] & {\c(a,x)} \arrow[lld, out=355,in=175,"\delta" description] \\
            & {\E(c,x)} \arrow[r]   & {\E(b,x)} \arrow[r] & {\E(a,x)} \arrow[lld,out=355,in=175]                       \\
            & {\E^2(a,x)} \arrow[r] & {}    \cdots              &                        
	\end{tikzcd}
\end{equation}
\begin{equation}\label{eq:leftles}
	\begin{tikzcd}[row sep=0.95em]
	{} & {\c(x,a)} \arrow[r]   & {\c(x,b)} \arrow[r] & {\c(x,c)} \arrow[lld,out=355,in=175,"\delta" description] \\
            & {\E(x,a)} \arrow[r]   & {\E(x,b)} \arrow[r] & {\E(x,c)} \arrow[lld,out=355,in=175] \\
            & {\E^2(x,a)} \arrow[r] & \cdots              &                      
\end{tikzcd}
\end{equation}
{The maps labeled by $\delta$ produce $\E$-extensions that are realized by good lifts.}
 
 \subsection{} \cite[Definition~3.7]{GNP23}
The extriangulated $\k$-category $(\c,\E)$ is \textit{0-Auslander} if it has enough projecitives, is hereditary, and any projective object $p\in\c$ is the domain of an inflation $\begin{tikzcd} p \arrow[r,tail] & q\end{tikzcd}$ with projective-injective codomain. 
A 0-Auslander extriangulated $\k$-category $(\c,\E)$ is \textit{reduced} if every projective-injective object therein is a zero object.

An exact dg $\k$-category $\A$ is called \textit{0-Auslander} if $(H_0\A,\E_{\A})$ is 0-Auslander as an extriangulated $\k$-category. A 0-Auslander exact dg $\k$-category $\A$ is \textit{reduced} if $(H_0\A,\E_{\A})$ is reduced.

\subsection{Example}\label{eg:0Aus} We include two classes of examples of 0-Auslander extriangulated $\k$-categories.
\begin{enumerate}
	\item\label{eg:0Aus_per2} \cite[§3.3.8]{GNP23} Let $A$ be a connective dg $\k$-algebra. Regarding $A$ as a dg $\k$-category with one object, let $H_0\per_{[0,1]}(A)$ be the two-term category of $A$, as defined in \nref{eg:exactdg}\ref{eg:exactdg1}. This is an example of a reduced 0-Auslander exact dg $\k$-category. 
	\item\label{eg:0Aus_mod} \cite[§3.3.2]{GNP23}
	For a positive integer $n$, let $\Lambda_n$ be the path $\k$-algebra of the linearly oriented $A_n$-quiver
	\begin{equation*}
	\begin{tikzcd}
	1\arrow[r] & 2 \arrow[r] & \cdots \arrow[r] & n.
	\end{tikzcd}
	\end{equation*}
	Then the abelian category $\fgMod(\Lambda_n)$ of finite-dimensional right $\Lambda_n$-modules is a 0-Auslander extriangulated $\k$-category.
	Among these, only $\fgMod(\Lambda_1)$ is reduced; for $n\geq 2$ there is exactly one non-zero projective-injective object in $\fgMod(\Lambda_n)$ up to isomorphism.
\end{enumerate}

 \subsection{}\label{def:silt}
 Let $(\c,\E)$ be an extriangulated $\k$-category.
A full subcategory $\calr\subseteq \c$ is called \textit{rigid (with respect to $\E$)} if it is closed under isomorphism and $\E^i(r,r')=0$ for all $i\geq 1$ and $r,r'\in \calr$. 
A rigid subcategory $\calr\subseteq \c$ is \textit{silting} if additionally $\thick(\calr)=\c$. 
Let $\rigid(\c)$ denote the set of rigid subcategories of $\c$, and $\silt(\c)$ the set of silting subcategories of $\c$. Given a rigid subcategory $\calr\subseteq \c$, let $\rigid_{\calr}(\c)$  (resp. $\silt_{\calr}(\c)$) denote subset of $\rigid(\c)$  (resp. $\silt(\c)$) consisting of the subcategories containing $\calr$.
An object $t\in \c$ is \textit{rigid} (resp. \textit{silting}) if the smallest additive subcategory of $\c$ containing $t$, denoted $\add(t)$, is {rigid} (resp. {silting}).

Let $\A$ be an exact dg $\k$-category. A full dg subcategory $\r\subseteq \A$ is called \textit{rigid} (resp. \textit{silting}) if $H_0\r$ is a rigid (resp. silting) subcategory of $H_0\A$. The notions of a \textit{rigid object} and \textit{silting object} in $\A$ are defined similarly.
In accordance with the notation introduced above, we define the sets $\rigid(\A)$, $\silt(\A)$, etc., in the same way. 

 \subsection{Definition} \label{def:g-finite}
Let $(\c,\E)$ be an extriangulated $\k$-category. We say that $(\c,\E)$ is \textit{$\bfg$-finite} if there are only finitely many silting subcategories of $\c$. An exact dg $\k$-category $\A$ is \textit{$\bfg$-finite} if $H_0\A$ is {$\bfg$-finite} as an extriangulated $\k$-category.
A connective dg $\k$-algebra $A$ is \textit{$\bfg$-finite}\footnote{For the special case of finite-dimensional $\k$-algebras, the original term is \textit{$\tau$-tilting finite} \cite[Definition 1.1]{DIJ19}. The term $\bfg$-finite arose from the study of the $\bfg$-vector fan.} if $\twoper(A)$ is $\bfg$-finite as a dg $\k$-category.

\subsection{}\label{AT22.5.12}
Let $(\c,\E)$ be an extriangulated $\k$-category.
Then the set $\silt(\c)$ admits the structure of a poset, where we set $\calr \geq \calr'$ if $\E^{i}(r,r')=0$ for all $i\geq 1$, $r\in\calr$ and $r'\in\calr'$ {\cite[Proposition 5.12]{AT22}}. In particular, if $\c$ is hereditary, we have that $\calr \geq \calr'$ precisely when $\E(r,r')=0$ for all $r\in\calr$ and $r'\in\calr'$.

\subsection{}\label{AT.5.4} If an extriangulated $\k$-category $(\c,\E)$ is Krull-Schmidt and contains a silting object $t$, then all silting subcategories of $\c$ are of the form $\add(s)$ for some silting object $s\in \c$ \cite[Proposition 5.4]{AT22}, and every silting object has the same number of indecomposable direct summands up to isomorphism \cite[Theorem 4.1]{AT23}. 

Suppose that the extriangulated $\k$-category $(\c,\E)$ is 0-Auslander. Then a rigid subcategory of $\c$ is silting precisely when it is maximal with respect to the rigidity property \cite[Theorem 4.3]{GNP23}.
Using the definition of projective and injective objects in \Cref{def:proj_inj}, it is easy to see that both $\proj(\c)$ and $\inj(\c)$ form maximally rigid subcategories of $\c$. It is also easy to see that $\proj(\c)$ is the maximal element in the poset $\silt(\c)$. Dually, we have that $\inj(\c)$ is the minimal element. When $\c$ is algebraic and enhanced by an exact dg $k$-category $\A$, we will denote the maximal and minimal elements in $\silt(\A)$ by $\p_{\A}$ and $\i_{\A}$, respectively.

\subsection{}\label{ff} 
{Let $\c$ be a $\k$-category. For an additive subcategory $\x\subseteq \c$, a \textit{right $\x$-approximation} of an object $y\in\c$ is a morphism $y^{\x}\to^{p}y$ in $\c$ such that $y^{\x}\in \x$ and the induced map $\c(x,y^{\x}) \to^{p\circ -}\c(x,y)$ surjects whenever $x\in \x$. We say that $\x$ is \textit{contravariantly finite} if every objects in $\c$ admits a right $\x$-approximation. One defines the notions of \textit{left $\x$-approximation} and \textit{covariant finiteness} dually.}

\subsection{}\label{Bongartz}
{Any rigid subcategory $\calr$ of a 0-Auslander extriangulated $\k$-category $(\c,\E)$ is a subcategory of a silting subcategory of $(\c,\E)$. One defines the \textit{maximal Bongartz completion} $\calr_{+}$ as the smallest additive subcategory of $\c$ containing $\calr$, the projective-injective objects in $(\c,\E)$, as well as all objects of the form $u_r^y$ appearing in a conflation}
\begin{equation*}
\begin{tikzcd}
	u_r^y \arrow[r,tail] & r' \arrow[r,"p_y",two heads] & y
\end{tikzcd}
\end{equation*}
{in which $y\in \inj(\c)$ and $p_y$ is a right $\add(\calr,\proj(\c)\cap \inj(\c))$-approximation. We say that $\calr$ \textit{admits a maximal Bongartz completion} if the deflation and right approximation $p_y$ exists for all injective objects $y\in \c$.
This is indeed a silting subcategory of $\c$ containing $\calr$ \cite[Corollary 4.8]{GNP23}, and it is the maximal such silting subcategory with respect to the partial order in \Cref{AT22.5.12} (this is shown using a long-exact sequence argument, see \eqref{eq:rightles}, and heredity). One defines the \textit{minimal Bongartz completion} $\calr_{\div}$ dually, and this is the minimal silting subcategory of $\c$ containing $\calr$. We say that $\calr$ \textit{admits Bongartz completions} if it admits a maximal and a minimal Bongartz completion.

\section{Silting reduction for extriangulated categories}\label{subsec:siltred}

 \subsection{Setup}\label{setup:siltred}
We fix an extriangulated $\k$-category $(\c,\E)$.
{For a full subcategory $\x$ of $\c$, let $\x^{\perp}$ (resp. ${^{\perp}\x}$) denote the full subcategory of $\c$ spanned by the objects $y\in \c$ with the property that $\E^{i}(x,y)=0$ (resp. $\E^{i}(y,x)=0$) for all $x\in\x$ and all positive integers $i$. We write $\x^{\perp_{\c}}$ (resp. ${^{\perp_{\c}}\x}$) when the choice of ambient extriangulated $\k$-category may be unclear.} Let $\x^{\vee}$ (resp. $\x^{\wedge}$) denote the smallest full subcategory of $\c$ containing $\x$ which is closed under extensions, cocones (resp. cones) and direct summands. 

Given a rigid subcategory $\calr\subseteq \c$, let $\z_{\calr}$ denote the full subcategory ${\calr^{\perp} \cap {^{\perp}\calr}} \subseteq \c$. One uses long-exact sequence arguments (\Cref{les}) to show that $\z_{\calr}$ is extension-closed in $\c$, which lets us induce an extriangulation on $\z_{\calr}$ using \Cref{lem:NP19.2.18}. 
The \textit{reduced category} of $(\c,\E)$ with respect to $\calr$ is given by ideal quotient 
$
	\overline{\c}_{\calr} \defeq \z_{\calr} / [\calr].
$
Since all objects in $\calr$ are projective-injective in $\z_{\calr}$,  the restriction of $\E$ to $\z_{\calr}$ induces a $\k$-linear bifunctor $\overline{\E}$ on $\overline{\c}_{\calr}$, yielding an extriangulated $\k$-category $(\overline{\c}_{\calr},\overline{\E})$. By \Cref{idquot}, the ideal quotient functor $\overline{(-)}\colon \z_{\calr} \to \overline{\c}_{\calr}$ becomes extriangulated when equipped with the natural homomorphism $\E \Longrightarrow \overline{\E}$ of functors $\overline{\c}_{\calr}\op\times \overline{\c}_{\calr}\to \Mod(\k)$.

\subsection{Proposition}\label{prop:siltred_0Aus}
Let $(\c,\E)$ be an extriangulated $\k$-category with enough projectives or injectives, and let $\calr\subseteq\c$ be a rigid subcategory. 
\begin{enumerate}
\item \label{prop:siltred_0Aus0} For any $i\geq 0$, the bifunctors $\overline{\E^{i}}$ and  $\overline{\E}^{i}$ are naturally isomorphic.
\item\label{prop:siltred_0Aus1} For all $x \in {^{\perp}\calr}$, $y\in \calr^{\perp}$, and $i >0$, the $\k$-vector space $\E^i(x,y)$ is trivial precisely when $\overline{\E}^i(x,y)$ is trivial.
\item\label{prop:siltred_0Aus2} The ideal quotient functor $\overline{(-)}\colon \z_{\calr} \to \overline{\c}_{\calr}$ induces a bijection
\begin{equation}\label{eq:rigidred_0Aus}
\overline{(-)}\colon \rigid_{\calr}(\c) \to \rigid(\overline{\c}_{\calr}).
\end{equation}
\item\label{prop:siltred_0Aus3}
Assume that $\silt_{\calr}(\z_{\calr}) = \silt_{\calr}(\c)$.
The map in \eqref{eq:rigidred_0Aus} then restricts to an isomorphism of posets (see \Cref{AT22.5.12})
\begin{equation}\label{eq:siltingred_0Aus}
\overline{(-)}\colon\silt_{\calr}(\c) \to \silt(\overline{\c}_{\calr}).
\end{equation}
\end{enumerate}
 \begin{proof}
 We prove \ref{prop:siltred_0Aus0} in the case where $(\c,\E)$ has enough projectives, noting that the case of enough injectives is dual. 
 We proceed by induction on $i$. The anchor step $i\leq 1$ trivially holds. Supposing that the claim holds for $i<j$, for some $j\geq 2$, we prove that the covariant functors $\overline{\E^{j}}(x,-)$ and  $\overline{\E}^{j}(x,-)$ are naturally isomorphic for any choice of $x\in\c$. Consider any conflation of the form
 \begin{equation}\label{eq:siltred_0Aus0_con}
 \begin{tikzcd}
x'\arrow[r,tail] & p \arrow[r,two heads] & x,
\end{tikzcd}
\end{equation}
where $p$ is projective, as in \eqref{eq:LN19.5.2}. It follows from \Cref{rem:her_highext} and the induction hypothesis that 
\begin{equation*}
	\overline{\E^{j}}(x,-) \cong \overline{\E^{j-1}}(x',-) \cong \overline{\E}^{j-1}(x',-) \cong \overline{\E}^{j}(x,-),
\end{equation*}
which establishes the inductive step and concludes the proof of \ref{prop:siltred_0Aus0}.
 
	
	We prove \ref{prop:siltred_0Aus1} under the assumptions that $(\c,\E)$ has enough projectives, the case of enough injectives being dual. We proceed by induction. To establish the anchor step $i=1$, it is sufficient to observe that the following are equiavalent:\footnote{Compare with Gorsky--Nakaoka--Palu \cite[Proposition 2.9(2)]{GNP23}. Although they work under stronger assumptions, namely $x,y\in \z_{\calr}$, the statement remains true and immediate from the definitions.}
	\begin{enumerate}
		\item $\E(x\oplus r_x,y\oplus r_y)=0$ for all $r_x,r_y\in \calr$,
		\item $\E(x,y)=0$,
		\item $\overline{\E}(x,y)=0$.
	\end{enumerate}
	For the inductive step, suppose that the assertion holds for all $i<j$, for some $j\geq 2$. Consider again the conflation in \eqref{eq:siltred_0Aus0_con}, where $x \in {^{\perp}\calr}$. Using that $p$ is projective, a long-exact sequence argument (see \eqref{eq:rightles}) shows that $x' \in {^{\perp}\calr}$. From \Cref{rem:her_highext} we get an isomorphism 
\begin{equation*}
	\overline{\E}^{j}(x,y) \cong \overline{\E}^{j-1}(x',y).
\end{equation*}
Thus, the vanishing of $\overline{\E}^{j}(x,y)$ is equivalent to the vanishing of $\overline{\E}^{j-1}(x',y)$.
By the induction hypothesis, the vanishing of $\overline{\E}^{j-1}(x',y)$ is equivalent to the vanishing of ${\E}^{j-1}(x',y)$. Using \Cref{rem:her_highext} again, we get another isomorphism 
\begin{equation*}
	{\E}^{j}(x,y) \cong {\E}^{j-1}(x',y),
\end{equation*}
whence the the vanishing of ${\E}^{j-1}(x',y)$ is equivalent to the vanishing of ${\E}^{j}(x,y)$. Putting it all together, the vanishing of ${\E}^{j}(x,y)$ is equivalent to the vanishing of $\overline{\E}^{j}(x,y)$, completing the inductive step. In conclusion, the vanishing of ${\E}^{i}(x,y)$ is equivalent to the vanishing of $\overline{\E}^{i}(x,y)$ for all $i\geq 1$.
		
	To prove \ref{prop:siltred_0Aus2}, we first note that the rigid subcategories of $\c$ containing $\calr$ are precisely the rigid subcategories of $\z_{\calr}$ containing $\calr$. This makes it possible for the ideal quotient functor $\overline{(-)}\colon \z_{\calr} \to \overline{\c}_{\calr}$ to induce a map out of $\rigid_{\calr}(\c)$.
	Since this ideal quotient functor induces a bijection from the set of full subcategories of $\z_{\calr}$ containing $\calr$ to the set of full subcategories of $\overline{\c}_{\calr}$, it now follows from \ref{prop:siltred_0Aus1} that the map in \eqref{eq:rigidred_0Aus} is a well-defined bijection.
	
	 To address \ref{prop:siltred_0Aus3}, we first show that a rigid subcategory of $\c$ containing $\calr$ is silting precisely when the corresponding rigid subcategory of $\overline{\c}_{\calr}$ is silting. 
	 It is straightforward to establish the equality 
	 \begin{equation}\label{eq:siltred_0Aus3_eq}
	 \overline{\thick_{\z_{\calr}}(\mathcal{X})} = \thick_{\overline{\c}_{\calr}}(\overline{\mathcal{X}})
	 \end{equation} 
	 of full subcategories of $\overline{\c}_{\calr}$, valid for any full subcategory $\mathcal{X}$ of $\z_{\calr}$.
	 Consequently, the bijection $$\overline{(-)}\colon \rigid(\z_{\calr}) \to \rigid_{\calr}(\overline{\c}_{\calr}),$$ 
	 namely that in \eqref{eq:rigidred_0Aus}, restricts to a bijection $\overline{(-)}\colon \silt_{\calr}(\z_{\calr}) \to \silt(\overline{\c}_{\calr})$. Having assumed that \\ $\silt_{\calr}(\z_{\calr}) = \silt_{\calr}(\c)$, we conclude that the map in \eqref{eq:siltingred_0Aus} indeed is a well-defined bijection. 
	 Furthermore, one now uses \ref{prop:siltred_0Aus1} to show that this map respects the partial order, rendering it an isomorphism of posets.
\end{proof}

\subsection{Definition} \cite[Definition 3.2]{Oga22} \cite[Definition 4.1]{NP19}\label{def:cotorpair}
Let $(\c,\E)$ be an extriangulated $\k$-category. Given two additive subcategories $\mathcal{U}$ and $\mathcal{V}$ of $\c$ that are closed under isomorphism and direct summands, consider the following conditions.
\begin{enumerate} 
\item[\textbf{(CT1)}]\label{item:CT1} For every $u\in\mathcal{U}$ and $v\in\mathcal{V}$, we have that $\E(u,v)=0$. 
\item[\textbf{(CT2)}]\label{item:CT2} Any morphism $u\to^{f}v$ in $\c$ with $u\in\mathcal{U}$ and $v\in\mathcal{V}$ factors through an object in $\mathcal{U}\cap \mathcal{V}$.
\item[\textbf{(CT3r)}]\label{item:CT3r} For every $x\in\c$, there exists a conflation in $\c$ of the form
\begin{equation}\label{eq:ct2r}
	\begin{tikzcd} x^{\mathcal{V}} \arrow[r,tail] & x^{\mathcal{U}} \arrow[r,two heads,"p"] & x, \end{tikzcd} 
\end{equation}
in which $p$ is a right $\mathcal{U}$-approximation and $x^{\mathcal{V}} \in \mathcal{V}$.
\item[\textbf{(CT3l)}]\label{item:CT3l} For every $x\in\c$, there exists a conflation in $\c$ of the form
\begin{equation}\label{eq:ct2l}
	\begin{tikzcd} x \arrow[r,tail,"i"] & x_{\mathcal{V}} \arrow[r,two heads] & x_{\mathcal{U}} \end{tikzcd} 
\end{equation}
in which $i$ is a left $\mathcal{V}$-approximation and $x_{\mathcal{U}}\in \mathcal{U}$.
\end{enumerate}
We say that the pair $(\mathcal{U},\mathcal{V})$ is a \textit{right cotorsion pair} (resp. \textit{left cotorsion pair}) in $(\c,\E)$ if \hyperref[item:CT2]{\textbf{(CT2)}} and \hyperref[item:CT3r]{\textbf{(CT3r)}} (resp. \hyperref[item:CT2]{\textbf{(CT2)}} and \hyperref[item:CT3l]{\textbf{(CT3l)}}) hold. The pair $(\mathcal{U},\mathcal{V})$ is a \textit{cotorsion pair} in $(\c,\E)$ \hyperref[item:CT1]{\textbf{(CT1)}}, \hyperref[item:CT3r]{\textbf{(CT3r)}}, and \hyperref[item:CT3l]{\textbf{(CT3l)}} hold.

\subsection{}\label{cotor_rems}
A pair which is both a left cotorsion pair and a right cotorsion pair need not be a cotorsion pair, since we are not requiring the condition \hyperref[item:CT1]{\textbf{(CT1)}} to be met. For example, if the extriangulated $\k$-category $\c$ non-semisimple (i.e. $\E$ is a non-zero bifunctor), then $(\c,\c)$ is a both a left and a right cotorsion pair, but not a cotorsion pair.

For a cotorsion pair $(\mathcal{U},\mathcal{V})$, it can be shown that $x\in \mathcal{U}$ if and only if $\E(x,v)=0$ for all $v\in \mathcal{V}$, and dually that $x\in \mathcal{V}$ if and only if $\E(u,x)=0$ for all $u\in \mathcal{U}$ \cite[Remark 4.4]{NP19}. 
Moreover, both $\mathcal{U}$ and $\mathcal{V}$ are extension-closed in $\c$ \cite[Remark 4.6]{NP19}. 

Note that \hyperref[item:CT2]{\textbf{(CT2)}} follows from \hyperref[item:CT3l]{\textbf{(CT3l)}} as soon as $\mathcal{U}$ is extension-closed, or dually from \hyperref[item:CT3r]{\textbf{(CT3r)}} as soon as $\mathcal{V}$ is extension-closed. Indeed, if $(\mathcal{U},\mathcal{V})$ is a cotorsion pair and $x \to^{f} v$ is a morphism in which $x\in \mathcal{U}$ and $v\in \mathcal{V}$, it factors through the object $x_{\mathcal{V}}$ in \eqref{eq:ct2l}. The intermediate object $x_{\mathcal{V}}$ is in ${\mathcal{V}}$ by definition and in ${\mathcal{U}}$ by extension-closure. In particular, the condition \hyperref[item:CT2]{\textbf{(CT2)}} is always met for cotorsion pairs.

\subsection{} \label{cotor_facts}
If $\mathcal{M}$ is a silting subcategory of $(\c,\E)$, then $({^{\perp}\mathcal{M}},\mathcal{M}^{\perp})$ is a cotorsion pair in $(\c,\E)$ such that ${^{\perp}\mathcal{M}} \cap \mathcal{M}^{\perp}= \mathcal{M}$ \cite[Proposition 5.10]{AT22}. 
Moreover, we have that $\mathcal{M}^{\vee} = {^{\perp} \mathcal{M}}$ and that $\mathcal{M}^{\wedge} = \mathcal{M}^{\perp}$. 
For a rigid subcategory $\calr$ of $\c$, we have that $\calr$ is silting in $\thick(\calr)$, whence we deduce that $\mathcal{R} = \mathcal{R}^{\vee} \cap \mathcal{R}^{\wedge}$. Futhermore, we have that $\mathcal{R}^{\perp} =  (\mathcal{R}^{\vee})^{\perp}$  and ${^{\perp} \mathcal{R}} = {^{\perp}(\mathcal{R}^{\wedge})}$ \cite[Lemma 4.4]{AT22}, from which one can derive that $\thick(\calr)\cap \calr^{\perp}=\calr^{\wedge}$ and that $\thick(\calr)\cap {^{\perp}\calr}=\calr^{\vee}$ \cite[Lemma 4.11]{AT22}. 
It follows that $\calr^{\vee}\cap \calr^{\perp} = \calr = {^{\perp}\calr}\cap \calr^{\wedge}$. One uses a long-exact sequence argument (see \Cref{les}) to show that $\calr^{\perp}$ is closed under cones, whence $\calr^{\wedge}\subseteq \calr^{\perp}$, and dually $\calr^{\vee}\subseteq {^{\perp}\calr}$ \cite[Lemma 4.3]{AT22}.

\subsection{Definition} \cite[Definition 3.2]{Oga22}
Let $(\c,\E)$ be an extriangulated $\k$-category.
An ordered tuple $\big((\mathcal{U},\mathcal{V})) , (\mathcal{X},\mathcal{Y}) \big)$ is called a \textit{generalized concentric twin cotorsion pair} (or \textit{gCTCP} for short)\footnote{Ogawa uses the terminology \textit{generalized Hovey twin cotorsion pair}, which is shortened to \textit{gHTCP}. We have altered the terminology for the following reason: the condition \textbf{(gCTCP1)} below encodes intertwining (hence \textit{twin}) and \textbf{(gCTCP2)} encodes the concentricity.} in $(\c,\E)$ if the pair $(\mathcal{U},\mathcal{V})$ is a left cotorsion pair in $(\c,\E)$, the pair $(\mathcal{X},\mathcal{Y})$ is a right cotorsion pair in $(\c,\E)$, and the following hold.
\begin{enumerate} 
\item[\textbf{(gCTCP1)}] $\mathcal{U} \subseteq \mathcal{X}$ and $\mathcal{V} \supseteq \mathcal{Y}$,
\item[\textbf{(gCTCP2)}] $\mathcal{U} \cap \mathcal{V}=\mathcal{X}\cap \mathcal{Y}$.
\end{enumerate}

\subsection{Condition \textbf{(gCP)}}\label{ass:cotor}
Let $(\c,\E)$ be an extriangulated $\k$-category. Consider the following condition on a rigid subcategory $\calr$ of $\c$: 
	\begin{enumerate}
		\item[]\label{ass:cotor0}  The ordered tuple $\big((\calr^{\vee}, \calr^{\perp}), ({^{\perp}\calr}, \calr^{\wedge})\big)$ is a gCTCP in $(\c,\E)$.\footnote{Since all subcategories involved are extension-closed, it follows from \Cref{cotor_rems} that \hyperref[item:CT2]{\textbf{(CT2)}} holds for both $(\calr^{\vee}, \calr^{\perp})$ and $({^{\perp}\calr}, \calr^{\wedge})$. By \Cref{cotor_facts}, it is thus equivalent to require that $(\calr^{\vee}, \calr^{\perp})$ satisfies \hyperref[item:CT3l]{\textbf{(CT3l)}} and that $({^{\perp}\calr}, \calr^{\wedge})$ satisfies \hyperref[item:CT3r]{\textbf{(CT3r)}}. Furthermore, both pairs happen to satisfy \hyperref[item:CT1]{\textbf{(CT1)}}.}
		\end{enumerate} 
	
\subsection{}\label{cotor_notech}
Although \nameref{ass:cotor} might seem technical, it applies to many interesting situations in practice. Consider first the case where $\c$ is a triangulated $\k$-category. 
One can reformulate Iyama--Yang's technical assumptions for silting reduction by imposing that both $(\calr^{\vee},\calr^{\perp})$ and $({^{\perp}\calr},\calr^{\wedge})$ be cotorsion pairs in $\c$ \cite[Condition \textbf{(CP)} and Remark 3.4]{LZZZ21} \cite[Proposition 3.2]{IY18}. 
By \Cref{cotor_facts}, this is a special instance of \nameref{ass:cotor}. 
Iyama--Yang point out that their technical assumption holds when $\c$ is Hom-finite, the rigid subcategory $\calr$ is of the form $\add(r)$ for some rigid object $r\in \c$, and there exists a silting object having $r$ as a direct summand. A particularly interesting example is $\c=H_0\per(A)$ for some connective dg $\k$-algebra $A$, and $\calr=\add(r)$ for some rigid object $r\in H_0\per(A)$.
In \nref{prop:0Aus_cotor} below, we prove that \nameref{ass:cotor} holds for all rigid subcategories of the two-term category $H_0\per_{[0,1]}(A)$. 

\subsection{Theorem}\label{prop:0Aus_cotor} 
Let $(\c,\E)$ be a 0-Auslander extriangulated $\k$-category and let $\calr$ be a rigid subcategory of $\c$. 
\begin{enumerate}
	\item\label{prop:0Aus_cotor_l}  The rigid subcategory $\calr$ admits a maximal Bongartz completion $\calr_{+}$ (see \Cref{Bongartz}) precisely when $(\calr, \calr^{\perp})$ is a left cotorsion pair.
	\item\label{prop:0Aus_cotor_r} Dually, the rigid subcategory $\calr$ admits a minimal Bongartz completion $\calr_{\div}$ precisely when $({^{\perp}\calr}, \calr)$ is a right cotorsion pair.
	\item\label{prop:0Aus_cotor_lr} The rigid subcategory $\calr$ admits Bongartz completions precisely when $\big((\calr, \calr^{\perp}),({^{\perp}\calr}, \calr)\big)$ is a gCTCP in $(\c,\E)$. In particular, \nameref{ass:cotorin} holds for all rigid subcategories of $\c$ admitting Bongartz completions.
\end{enumerate} 
\begin{proof}
	It suffices to prove the assertion in \ref{prop:0Aus_cotor_l}, since \ref{prop:0Aus_cotor_r} is dual, and \ref{prop:0Aus_cotor_lr} easily follows from \ref{prop:0Aus_cotor_l} and \ref{prop:0Aus_cotor_r}. Since $\calr$ is rigid, we have that $\calr\subseteq \calr^{\perp}$, whence \hyperref[item:CT2]{\textbf{(CT2)}} holds for the pair $(\calr, \calr^{\perp})$. Thus, it suffices to show that $\calr$ admits a maximal Bongartz completion precisely when \hyperref[item:CT3l]{\textbf{(CT3l)}} holds for the pair $(\calr, \calr^{\perp})$.
	
	Suppose that $\calr$ admits a maximal Bongartz completion $\calr_{+}$. 
	Letting $x\in\c$ be arbitrary, one uses the cotorsion pair $({^{\perp}\calr_{+}},\calr_{+}^{\perp})$ to construct a conflation
	\begin{equation*}
	\begin{tikzcd}
		x\arrow[r,tail] & x_{(\calr_{+}^{\perp})} \arrow[r,two heads] & x_{(^{\perp}\calr_{+})}
	\end{tikzcd}
	\end{equation*}
	where $ x_{(\calr_{+}^{\perp})} \in {\calr_{+}^{\perp}} $ and $ x_{({^{\perp}\calr_{+}})} \in {^{\perp}\calr_{+}}$.
	It follows from heredity and a long-exact sequence argument (along the lines of \eqref{eq:leftles}) that $x_{({^{\perp}\calr_{+}})} \in {{\calr_{+}^{\perp}}}$, whence $x_{({^{\perp}\calr_{+}})} \in  {^{\perp}\calr_{+}} \cap {{\calr_{+}^{\perp}}}  =  \calr_{+}$ by \Cref{cotor_facts} and the fact that $\calr_{+}$ is a silting subcategory of $\c$.
	By the definition of the maximal Bongartz completion $\calr^{+}$, one can thus pick a conflation
	\begin{equation}\label{eq:bonleft}
	\begin{tikzcd}
		{x_{({^{\perp}\calr_{+}})}} \arrow[r,tail,"\iota'"] & x'  \arrow[r,two heads] & y,
	\end{tikzcd}
	\end{equation}
	in which $y$ is an injective object in $\c$ and $x' \simeq x'_{(\calr)}\oplus n$ is a biproduct of an object $x'_{(\calr)}$ in $\calr$ with a projecitive-injective object $n$ in $\c$. By heredity again, the map $\E(x',x) \to^{\E(\iota',x)} \E(x_{({^{\perp}\calr_{+}})},x)$ is an epimorphism of $\k$-modules, whence there is a diagram of conflations
	\begin{equation*}
	\begin{tikzcd}
x \arrow[r, tail] \arrow[d, equal] & x_{(\calr_{+}^{\perp})} \arrow[r, two heads] \arrow[d, tail] \arrow[rd, phantom] & x_{({^{\perp}\calr_{+}})} \arrow[d, tail,"\iota'"] \\
x \arrow[r, tail]                    & z \arrow[r,two heads] \arrow[d, two heads]                                                             & x' \arrow[d, two heads]                \\
                                     & y \arrow[r, equal]                                                                         & y                                     
\end{tikzcd}
\end{equation*}
where the conflation in \eqref{eq:bonleft} is shown along the rightmost column, and the two upper rows is a good lift of the middle row along $\iota'$ (indeed, the map $\delta$ in \eqref{eq:rightles} produces good lifts). Since $z$ is an extension of $y\in \c^{\perp}$ in $x_{(\calr_{+}^{\perp})}\in \calr_{+}^{\perp}\subseteq \calr^{\perp}$, we have $z\in \calr^{\perp}$ by extension-closure. Further, since the projective-injective object $n$ is a direct summand of $z$, we can remove $n$ from the middle and right terms of the middle row, yielding a conflation
\begin{equation}\label{eq:prop:0Aus_cotor_final}
\begin{tikzcd}
	x \arrow[r, tail,"\iota''"]                    & z/n \arrow[r,two heads]                                                              & x'_{(\calr)}    
\end{tikzcd}
\end{equation}
where $z/n\in \calr^{\perp}$ and $x'_{(\calr)}\in \calr$. The morphism $\iota''$ is a left $\calr^{\perp}$-approximation since $x'_{(\calr)}\in \calr$, as is shown by a long-exact sequence argument (see \eqref{eq:rightles}). 
We may now conclude, since the existence of the conflation in \eqref{eq:prop:0Aus_cotor_final} proves that \hyperref[item:CT3l]{\textbf{(CT3l)}} holds for $(\calr, \calr^{\perp})$.

Conversely, suppose that \hyperref[item:CT3l]{\textbf{(CT3l)}} holds for the pair $(\calr, \calr^{\perp})$, and let $y\in \c$ be an injective object. Since the definition of 0-Auslander extriangulated $\k$-categories is self-dual, there exists a deflation $\begin{tikzcd} q\arrow[r,two heads,"\pi"] & y,\end{tikzcd}$ where $q$ is projective-injective \cite[Proposition 3.6]{GNP23}. One consequently produces a conflation
\begin{equation*}
\begin{tikzcd}
p \arrow[r,"\iota",tail] & q \arrow[r,"\pi", two heads]  & y.
\end{tikzcd}
\end{equation*}
Consider the diagram of conflations
\begin{equation*}
	\begin{tikzcd}
		p \arrow[d, tail,"\iota"] \arrow[r, tail]      & p_{(\calr^{\perp})} \arrow[r, two heads] \arrow[d, tail] & p_{(\calr)} \arrow[d, equal] \\
		q \arrow[d, two heads,"\pi"] \arrow[r, tail] & b \arrow[r, two heads] \arrow[d, two heads,"\pi' "]              & p_{(\calr)}                    \\
		y \arrow[r, equal]                   & y                                                        &                               
	\end{tikzcd}
\end{equation*}
where we have used \hyperref[item:CT3l]{\textbf{(CT3l)}} to produce the top row, and the middle row is a good lift of the top row along $\iota$. The projective-injectivity of $q$ renders the middle row a split conflation, so $b\in \add(\calr,\proj(\c)\cap\inj(\c))$. Since $p_{(\calr^{\perp})}\in \calr^{\perp} = \add(\calr,\proj(\c)\cap\inj(\c))^{\perp}$, one uses a long-exact sequence argument (see \eqref{eq:leftles}) to show that $\pi'$ is a right $ \add(\calr,\proj(\c)\cap\inj(\c))$-approximation. We have shown that any injective object in $\c$ admits a right $\add(\calr,\proj(\c)\cap\inj(\c))$-approximation which is a deflation, which is precisely what is needed for $\calr$ to admit a maximal Bongartz completion.
\end{proof} 

\subsection{Proposition}\label{lem:red_0Aus}
Let $(\c.\E)$ be an extriangulated $\k$-category, let $\calr\subseteq\c$ be a rigid subcategory, and let $\z_{\calr}$ and $\overline{\c}_{\calr}$ be as in \nref{setup:siltred}. Then:\footnote{There exist results in the literature proving \ref{lem:red_0Aus1.5} and \ref{lem:red_0Aus2} below in the case where $(\c,\E)$ is a two-term category and $\calr$ is either the subcategory of projectives or the subcategory of injectives \cite[§4]{FGPPP23}.}
\begin{enumerate}
	\item\label{lem:red_0Aus1} If $(\c,\E)$ is hereditary, so are $\z_{\calr}$ and $\overline{\c}_{\calr}$.
	\item\label{lem:red_0Aus1.25} Suppose that $(\c,\E)$ is 0-Auslander and admits Bongartz completions. Then  $\calr_{+}^{\perp}=\calr^{\perp}$ and ${^{\perp}\calr_{\div}}={^{\perp}\calr}$. 
	\end{enumerate}
	Assume henceforth that $(\c,\E)$ is as in \ref{lem:red_0Aus1.25}.
	\begin{enumerate}
	\setcounter{enumi}{2}
	\item\label{lem:red_0Aus1.5} The projective objects in $\z_{\calr}$ are precisely those in $\calr_{+}$, and the injective objects in $\z_{\calr}$ are precisely those in $\calr_{\div}$.
	\item\label{lem:red_0Aus2} The extriangulated categories $\z_{\calr}$ and $\overline{\c}_{\calr}$ are 0-Auslander. 
	\item\label{lem:red_0Aus3} If the rigid subcategory $\calq\subseteq \c$ containing $\calr$ admits Bongartz completions, then the rigid subcategory $\overline{\calq}$ of $ \overline{\c}_{\calr}$ admits Bongartz completions.
\end{enumerate}
\begin{proof}
	The statement in \ref{lem:red_0Aus1} is known to experts, but we give a brief account of the proof here. Since $\z_{\calr}$ is extension-closed in the hereditary extriangulated $\k$-category $(\c,\E)$, it is straightforward to show that the restriction of $\E$ forms a hereditary extriangulated $\k$-category $(\z_{\calr},\E)$ \cite[Lemma 2.4]{GNP23}. The heredity of $\overline{\c}_{\calr}$ is now easily deduced using the definitions \cite[Proposition 2.9(1)]{GNP23}.

	We prove the former assertion of \ref{lem:red_0Aus1.25}, the latter being dual. Since $\calr$ is a subcategory of $\calr_{+}$, it is immediate that $\calr_{+}^{\perp}\subseteq\calr^{\perp}$. To prove the converse, we fix $x\in \calr^{\perp}$ and $u\in \calr_{+}$. By the definition of the maximal Bongartz completion in \Cref{Bongartz}, there exists a conflation
	 \begin{equation*}
	\begin{tikzcd}
	u\arrow[r,tail] & q \arrow[r,two heads] & y
	\end{tikzcd}
	\end{equation*}
where $y$ is injective in $\c$ and $q = r \oplus n$ is a biproduct of an object $r\in\calr$ and a projective-injective object $n\in \c$. By \Cref{les}, one induces a long-exact sequence
	\begin{equation*}
	\begin{tikzcd}
	\cdots & \arrow[l] \E^2(i,x) \arrow[l] & \E(u,x)\arrow[l] & \E(q,x) \arrow[l] &  \arrow[l] \cdots
	\end{tikzcd}
	\end{equation*}
	The rightmost term shown vanishes by assumption, and the leftmost term vanishes by heredity. The middle term thus vanishes. Since $u$ was arbitrarily chosen, we have shown that $x\in \calr_{+}^{\perp}$, as desired.

	We now prove the assertion in \ref{lem:red_0Aus1.5} concerning $\calr_{+}$. Fixing a projective object $p$ in $\z_{\calr}$, consider the following conflation in $\c$, given by the cotorsion pair $({^{\perp}\calr_{+}},\calr_{+}^{\perp})$.
	\begin{equation}\label{eq:equiv_H0AP0_cotor_conf}
	\begin{tikzcd}
	p^{(\calr_{+}^{\perp})} \arrow[r,tail] & p^{({^{\perp}\calr_{+}})} \arrow[r,two heads] & p.
	\end{tikzcd}
\end{equation}
Since $\c$ is hereditary and $p^{({^{\perp}\calr_{+}})}\in {^{\perp}\calr_{+}} \subseteq {^{\perp}\calr}$, one uses \Cref{les} to show that $p^{(\calr_{+}^{\perp})}\in {^{\perp}\calr}$, implying that $$p^{(\calr_{+}^{\perp})}\in {^{\perp}\calr}\cap {\calr_{+}^{\perp}} \subseteq {^{\perp}\calr}\cap {\calr_{}^{\perp}} = \z_{\calr}.$$ By the extension-closure of $\z_{\calr}$, the entire conflation in \eqref{eq:equiv_H0AP0_cotor_conf} therefore lies in $\z_{\calr}$. Moreover, \begin{equation*}
\begin{tikzcd}[column sep=1em]
p^{({^{\perp}\calr_{+}})}  \in {^{\perp}\calr_{+}} \cap \z_{\calr} \arrow[r,symbol=\subseteq]& {^{\perp}\calr_{+}} \cap \calr^{\perp} \arrow[r,equal,"\ref{lem:red_0Aus1.25}"'] &{^{\perp}\calr_{+}} \cap \calr_{+}^{\perp} \arrow[r,equal,"\Cref{cotor_facts}"']& \calr_{+}.
\end{tikzcd}
\end{equation*}
The projectivity of $p$ in $\z_{\calr}$ renders the conflation in  \eqref{eq:equiv_H0AP0_cotor_conf} split, whence the object $p$ is a direct summand of $p^{({^{\perp}\calr_{+}})}$, placing $p$ in $\calr_{+}$. Since $p$ was arbitrarily chosen, we have shown that all projective objects in $\z_{\calr}$ are in $\calr_{+}$. Conversely, all objects in $\calr_{+}$ are projective in $\z_{\calr}$, since $\calr_{+}^{\perp_{\z_{\calr}}} = \calr_{+}^{\perp_{\c}}\cap \z_{\calr} = \calr^{\perp_{\c}}\cap \z_{\calr} = \z_{\calr}$, the second equality following from \ref{lem:red_0Aus1.25} above. 

	We move on to prove \ref{lem:red_0Aus2}. Firstly, we show that $\z_{\calr}$ is 0-Auslander.
	By \ref{lem:red_0Aus1} above, it suffices to prove that every projective object $p\in \z_{\calr}$ appears in a conflation in $\z_{\calr}$
\begin{equation}\label{eq:equiv_H0AP0_cotor_conf'}
	\begin{tikzcd}
	p \arrow[r,tail] & q \arrow[r,two heads] & y,
	\end{tikzcd}
\end{equation}
where $q$ is projective-injective. 
By \ref{lem:red_0Aus1.5}, we have that $p \in \calr_{+}$. The definition of the maximal Bongartz completion in \Cref{Bongartz} gives a conflation in $\c$
\begin{equation}\label{eq:equiv_H0AP0_cotor_conf''}
	\begin{tikzcd}
	p \arrow[r,tail] & q' \arrow[r,two heads] & j,
	\end{tikzcd}
\end{equation}
where $j$ is injective in $\c$ and $q' = r \oplus n$ is a biproduct of an object $r\in\calr$ and a projective-injective object $n\in \c$. Now, \nref{prop:0Aus_cotor}\ref{prop:0Aus_cotor_r} states that $({^{\perp}\calr},\calr)$ is a right cotorsion pair in $(\c,\E)$, letting us form the following diagram of conflations
\begin{equation*}
	\begin{tikzcd}
                                     & j^{(\calr)} \arrow[d, tail] \arrow[r, equal] & j^{(\calr)} \arrow[d, tail]       \\
	p \arrow[r, tail] \arrow[d, equal] & q \arrow[d, two heads] \arrow[r, two heads]             & j^{({^{\perp}\calr})} \arrow[d, two heads] \\
	p \arrow[r,tail]                          & q' \arrow[r, two heads]                            & j                                         
	\end{tikzcd}
\end{equation*}
in which the bottom row is the conflation in \eqref{eq:equiv_H0AP0_cotor_conf''}, and the middle row and middle column are good lifts. Since $j\in \c^{\perp}\subseteq \calr^{\perp}$ and $j^{(\calr)}\in \calr \subseteq \calr^{\perp}$, extension-closure gives that $j^{({^{\perp}\calr})} \in \calr^{\perp}$, whence $j^{({^{\perp}\calr})} \in \z_{\calr}$. Consequently, the middle row is a conflation in $\z_{\calr}$. Writing $q \simeq t \oplus n$, where $n$ is the projective-injective direct summand of $q'=r\oplus n$ fixed above, we have that $t$ is in $\calr$, being an extension of $r\in \calr$ in $j^{(\calr)}\in\calr$ vertically.
We deduce that the middle row is the sought conflation in \eqref{eq:equiv_H0AP0_cotor_conf'}, and that $\z_{\calr}$ is 0-Auslander. 
Furthermore, the ideal $[\calr]$ of $\z_{\calr}$ is, by definition, generated by morphisms with injective domain and projective codomain. Since ideal quotients of this form preserve the 0-Auslander property \cite[Corollary 3.3]{FGPPP23}, it follows that the ideal quotient $\overline{\c}_{\calr} = \z_{\calr}/[\calr]$ indeed is 0-Auslander whenever $(\c,\E)$ is. 

Finally, we prove \ref{lem:red_0Aus3}. We only prove that the maximal Bongartz completion of $\overline{\calq}$ in $\overline{\c}_{\calr}$ exists. Since $\calq$ is a subcategory of the 0-Auslander extriangulated $\k$-category $\z_{\calr}$, it is easy to see that the following suffices: the rigid subcategory $\calq\subseteq \z_{\calr}$ admits a maximal Bongartz completion. 
By \ref{lem:red_0Aus1.25} above, any injective object $y$ in $\z_{\calr}$ belongs to the minimal Bongartz completion $\calr_{\div}$ of $\calr$ in $\c$. 
Consider the following commutative diagram of conflations
 \begin{equation*}
\begin{tikzcd}
p \arrow[r, "\alpha"] \arrow[d]                               & r' \arrow[r, two heads] \arrow[d]           & y \arrow[d, equal] \\
p_{(\calq^{\perp})} \arrow[d, two heads] \arrow[r] & q \arrow[d, two heads] \arrow[r, two heads,"\beta"] & y                    \\
p_{(\calq)} \arrow[r, equal]                                & p_{(\calq)}                                 &                     
\end{tikzcd}
\end{equation*}
in which $p$ is projective in $\c$, the morphism $\alpha$ is a left $\add(\calr, \mathcal{P}\cup\mathcal{I})$-approximation (see \Cref{Bongartz}), the left column comes from the left cotorsion pair $({\calq},\calq^{\perp})$ in $(\c,\E)$ (see \nref{prop:0Aus_cotor}), and the rest of the diagram is obtained from good lifts. Since $p_{({\calq^{\perp}})}$ is an extension of $p_{(\calq)}\in \calq \subseteq {^{\perp}\calr}$ in $p\in {^{\perp}\c}\subseteq {^{\perp}\calr}$, it is $\calq^{\perp}\cap {^{\perp}\calr}$, and thus also in $\z_{\calr}$. The object $q$ is in $\add(\calq, \mathcal{P}\cup\mathcal{I})$, since it is an extension of $p_{(\calq)}\in \calq$ in $r' \in \add(\calr, \mathcal{P}\cup\mathcal{I})\subseteq  \add(\calq, \mathcal{P}\cup\mathcal{I})$. 
Since $\E(x,p_{(\calq)})=0$ for all $x\in \add(\calq, \mathcal{P}\cup\mathcal{I})$, it now follows from a long-exact sequence argument (see \eqref{eq:leftles}) that $\beta$ a right $\add(\calq, \mathcal{P}\cup\mathcal{I})$-approximation, whence the middle row gives the means for constructing the maximal Bongartz completion $\calq_{+}$ in $\z_{\calr}$. We conclude that $\overline{\calq}$ admits a maximal Bongartz completion in $\overline{\c}_{\calr}$.
\end{proof}

\subsection{}\label{redcat}
Let $(\c,\E)$ be an extriangulated $\k$-category with enough projectives or injectives, and let $\mathcal{N}\subseteq \c$ denote the subcategory of projective-injective objects. Note that $\z_{\mathcal{N}}=\c$.
One defines the \textit{reduced extriangulated $\k$-category} of $(\c,\E)$ by $(\overline{\c}_{\mathcal{N}},\overline{\E}) = (\c/[\mathcal{N}],\overline{\E})$. This is clearly a reduced extriangulated $\k$-category. If $(\c,\E)$ is 0-Auslander, so is $(\overline{\c}_{\mathcal{N}},\overline{\E})$, by \nref{lem:red_0Aus}\ref{lem:red_0Aus2}. 
Moreover, all silting subcategories of the 0-Auslander extriangulated $\k$-category $\c$ contain $\mathcal{N}$ \cite[Remark 4.4(d)]{GNP23}, whence the ideal quotient functor $\overline{(-)}\colon \c \to \overline{\c}_{\mathcal{N}}$ induces a poset isomorphism
\begin{equation*}
	\silt(\c) = \silt_{\mathcal{N}}(\c) \to^{\sim} \silt(\overline{\c}_{\mathcal{N}}).
\end{equation*}
This poset isomorphism has previously been constructed by Pan--Zhu \cite[Theorem 1.2 and Lemma 4.14]{PZ24}.

\subsection{Lemma}\label{cor:0AusZ} Let $(\c,\E)$ be an extriangulated $\k$-category with enough projectives or injectives. For a rigid subcategory $\calr\subseteq \c$ satisfying \nameref{ass:cotor}, we have that $\silt_{\calr}(\z_{\calr}) = \silt_{\calr}(\c)$. In particular, the assertion in \nref{prop:siltred_0Aus}\ref{prop:siltred_0Aus3} holds.
\begin{proof}
	We refer to Liu--Zhou--Zhou--Zhu \cite[Lemma 3.11]{LZZZ21}. Although their condition \textbf{(CP)} is stronger than our \nameref{ass:cotor}, the reader may check that \nameref{ass:cotor} is sufficient to conduct this proof.\footnote{For the special case where $(\c,\E)$ is 0-Auslander and admits Bongartz completions, one can alternatively proceed as follows: It was asserted in \Cref{AT.5.4} that a rigid subcategory of $\c$ is silting precisely when it is maximal with respect to the rigidity property. By \nref{lem:red_0Aus}\ref{lem:red_0Aus2}, the same can be said of rigid subcategories of $\z_{\calr}$. Recalling that $\rigid_{\calr}(\z_{\calr}) = \rigid_{\calr}(\c)$, we now conclude, since a rigid subcategory of $\c$ is maximal if and only if it is maximal as a rigid subcategory of $\z_{\calr}$.}
\end{proof}

\subsection{Theorem}\label{cor:siltred0Aus}
Let $(\c,\E)$ be an extriangulated $\k$-category, and let $\calr\subseteq \c$ be a rigid subcategory.
\begin{enumerate}
	\item\label{prop:siltred} Suppose that $(\c,\E)$ has enough projectives or injectives, and let $\z_{\calr}$ and $\overline{\c}_{\calr}$ be as in \nref{setup:siltred}. Assume that \nameref{ass:cotor} holds for $\calr$.
Then the composite functor
	\begin{equation*}
	\begin{tikzcd}
		{\z_{\calr}} \arrow[r,hook] & \c \arrow[r,"L_{\calr}"] & \c/\thick(\calr),
	\end{tikzcd}
	\end{equation*}
	where $L_{\calr}$ is the Verdier localization (see \Cref{NOSloc}), induces an equivalence $\overline{\c}_{\calr} \to \c/\thick(\calr)$. Consequently, the $\k$-category $\c/\thick(\calr)$ can be equipped with an extriangulation such that $L_{\calr}$ carries the structure of extriangulated functor.
	\item\label{cor:siltred0Aus1}  If $(\c,\E)$ and $\calr$ are as in \ref{prop:siltred}, then Verdier localization functor $ \c \to^{L_{\calr}} \c/\thick(\calr)$ induces a bijection
\begin{equation*}
\rigid_{\calr}(\c) \to \rigid(\c/\thick(\calr))
\end{equation*}
which restricts to an isomorphism of posets
\begin{equation*}
\silt_{\calr}(\c) \to \silt(\c/\thick(\calr)).
\end{equation*}
	\item\label{cor:siltred0Aus2} If $(\c,\E)$ is 0-Auslander and admits Bongartz completions, then the Verdier quotient $\c/\thick(\calr)$ is 0-Auslander and admits Bongartz completions.  	
\end{enumerate}
\begin{proof}
	After we have proved \ref{prop:siltred}, the assertions in \ref{cor:siltred0Aus1} and \ref{cor:siltred0Aus2} follow from \nref{cor:0AusZ}, \nref{prop:siltred_0Aus}, \nref{prop:0Aus_cotor}, \nref{lem:red_0Aus}\ref{lem:red_0Aus2}, and \nref{lem:red_0Aus}\ref{lem:red_0Aus3}.
	Since $L_{\calr}$ sends every object in $\calr$ to a zero object, we can indeed use induce a functor $\Phi\colon\overline{\c}_{\calr} \to \c/\thick(\calr)$. To prove our result, we show that this functor admits a quasi-inverse, and that the Verdier localization $L_{\calr}$ can be factored as $\c\to^{F} \overline{\c}_{\calr}\to^{\Phi} \c/\thick(\calr)$ up to natural isomorphism, where $F$ is an extriangulated functor.
	The result essentially follows from a work of Ogawa \cite[Theorem 3.10]{Oga22}. Below, we spell out Ogawa's approach, in order to convince the reader that all necessary assumptions of are met for the situation at hand.
	
	We have assumed that $\big((\calr^{\vee}, \calr^{\perp}), ({^{\perp}\calr}, \calr^{\wedge})\big)$ is a gCTCP in $(\c,\E)$. 
	The left cotorsion pair $(\calr^{\vee}, \calr^{\perp})$ determines a functor $\lambda\colon \c/[\calr] \to \calr^{\perp}/[\calr]$, where $\lambda(x) = x_{(\calr^{\perp})}$ (see \hyperref[item:CT3l]{\textbf{(CT3l)}}), which is a left adjoint of the inclusion of $\calr^{\perp}/[\calr]$ into $\c/[\calr]$ \cite[Lemma 3.3]{Oga22}. Dually the right cotorsion pair $({^{\perp}\calr}, \calr^{\wedge})$ determines a functor $\rho\colon \c/[\calr] \to {^{\perp}\calr}/[\calr]$, where $\rho(x) = x^{({^{\perp}\calr})}$ (see \hyperref[item:CT3r]{\textbf{(CT3r)}}), which is a right adjoint of the inclusion of ${^{\perp}\calr}/[\calr]$ into $\c/[\calr]$. 
	Consider the composite functors $\rho\lambda$ and $\lambda\rho$.
	Letting $\overline{\rho\lambda}$ and $\overline{\lambda\rho}$ denote the composites of these functors with the ideal quotient functor $\overline{(-)}\colon \c\to {\c}/{[\calr]}$, the first step of our proof will be to show that $\overline{\rho\lambda}$ and $\overline{\lambda\rho}$ are naturally isomorphic, and that their image is in $\overline{\c}_{\calr}$. For an arbitrary $x\in \c$, consider the solid part of the diagram
	\begin{equation}\label{eq:diadia}
		\begin{tikzcd}
x^{(\calr^{\wedge})} \arrow[r, tail] \arrow[d, dashed,"1"] & x^{({^{\perp}\calr})} \arrow[r, two heads] \arrow[d, dashed, tail] & x \arrow[d, tail]                             \\
x^{(\calr^{\wedge})} \arrow[r, dashed, tail]           & z \arrow[d, two heads, dashed] \arrow[r, two heads, dashed]      & x_{(\calr^{\perp})} \arrow[d, two heads] \\
                                                      & x_{(\calr^{\vee})} \arrow[r, dashed,"1"]                            & x_{(\calr^{\vee})}                        
\end{tikzcd}
	\end{equation}
	Since $\calr$ is silting in $\thick(\calr)$, it follows form \Cref{cotor_facts} that $\E^2(x_{(\calr^{\vee})},x^{(\calr^{\wedge})})=0$, so we can complete the diagram in such a way that top two rows give a good lift of the middle row, and the arrows labeled by $1$ are identity morphisms. It is also the case that the two rightmost columns give a good lift of the middle column. The object $z$ is an extension of $x_{(\calr^{\vee})}$ in $x^{({^{\perp}\calr})}$ vertically, placing it in ${^{\perp}\calr}$, since ${^{\perp}\calr}$ is extension-closed and contains $\calr^{\vee}$. Working horizontally, one shows that $z\in \calr^{\perp}$, whence $z\in \z_{\calr}$. 
	Moreover, we can express  $z$ as either ${x_{(\calr^{\perp})}}^{(\calr^{\perp})}= \overline{\rho\lambda}(x)$ or ${x^{({^{\perp}\calr})}}_{(\calr^{\perp})}=\overline{\lambda\rho}(x)$. It is now straightforward to deduce that $\overline{\rho\lambda}$ and $\overline{\lambda\rho}$ are naturally isomorphic as functors from $\c$ to $\overline{\c}_{\calr}$. Henceforth, we denote $\overline{\rho\lambda}$ by $F$.
	
	We will now show that the functor $\lambda\colon \c/[\calr]\to {\calr^{\perp}}/[\calr]$ is extriangulated. A dual argument will show that $\rho$ is extriangulated, and consequently $F\colon \c\to\overline{\c}_{\calr}$ will be extriangulated, since the ideal quotient functor $\overline{(-)}\colon \c \to \c/[\calr]$ is extriangulated by \Cref{idquot}. For an arbitrary $x\in \c$, the left cotorsion pair $(\calr^{\vee},\calr^{\perp})$ gives a conflation
	\begin{equation}\label{confyconf}
	\begin{tikzcd}
		x \arrow[r,tail,"\eta_x"]  & x_{(\calr^{\perp})} \arrow[r, two heads] & x_{(\calr^{\vee})},
	\end{tikzcd}	
	\end{equation}	
	in which the inflation $\eta_x$ is a left ${\calr^{\perp}}$-approximation.
	The morphism $\overline{\eta_x}$ in $\c/[\calr]$ is then the unit of the adjunction $\lambda \dashv i$, where $i$ is the inclusion functor of $\calr^{\perp}/[\calr]$ into $\c/[\calr]$.  
	By \Cref{cotor_facts} and \Cref{prop:siltred_0Aus}\ref{prop:siltred_0Aus1}, we have that $\overline{\E}^i(x_{(\calr^{\vee})},y)=0$, for all $x\in \c$, all $y\in \calr^{\perp}$, and $i\in \{1,2\}$, whence
	$\begin{tikzcd}[column sep=2.5em]
	\overline{\E}(\lambda-,-) \arrow["{\overline{\E}(\eta,-)}",r] &\overline{\E}(-,-)
\end{tikzcd}$ is a natural isomorphism of functors from ${\c \over [\calr]} \times {{\calr^{\perp}}\over[\calr]}$ to $\Mod(\k)$. Consider now the natural transformation $\psi$ given by the composite
	\begin{equation*}
\begin{tikzcd}[column sep=4em]
	\overline{\E}(-,-) \arrow["{\overline{\E}(-,\eta-)}",r] & \overline{\E}(-,\lambda-) \arrow[r,"{\overline{\E}(\eta-,-)^{-1}}"] & \overline{\E}(\lambda-,\lambda-).
\end{tikzcd}
\end{equation*}
We will show that $(\lambda,\psi)$ is an extriangulated functor. Let $\xi\in \overline{\E}(c,a)$ be an $\overline{\E}$-extension realized by the top row in the diagram in \eqref{lem:extadjcom} below.
\begin{equation}\label{lem:extadjcom}
	\begin{tikzcd}
	a \arrow[r,tail] \arrow[d, "\eta_a"] & b \arrow[r, two heads] \arrow[d,"\kappa_1"]  & c \arrow[d, equal]            \\
	\lambda a \arrow[r,tail]                    & b' \arrow[r, two heads] \arrow[d,"\kappa_2"]          & c   \arrow[d, "\eta_c"]                            \\
	\lambda a \arrow[u, equal] \arrow[r,tail]  & b'' \arrow[r, two heads]& \lambda c 
	\end{tikzcd}
\end{equation}
The rest of the diagram is constructed by first choosing a good lift along $\eta_a$, and then using that ${\E(\eta,-)}$ is a natural isomorphism to add the bottom row. It is to be shown that the square in \eqref{eq:extrifun} commutes, which amounts to that the bottom row is isomorphic to $\lambda\xi$ in $\mathcal{H}^{}_{3\mathrm{t}}(\calr^{\perp}/[\calr])$. It suffices to show that the morphism $\kappa_2\kappa_1$ is a left $\calr^{\perp}$-approximation. The object $b''$ is in $\calr^{\perp}$ thanks to extension closure along the bottom row in \eqref{lem:extadjcom}. For an object $y\in \calr^{\perp}$ we use \Cref{les} to induce a morphism of long-exact sequences
\begin{equation*}
	\begin{tikzcd}
\cdots & {\E(c,y)} \arrow[l]                                     & {\c(a,y)} \arrow[l]                                                & {\c(b,y)} \arrow[l]                                         & {\c(c,y)} \arrow[l]                                                \\
\cdots & {\E(\lambda c,y)} \arrow[l] \arrow[u, "{\E(\eta_c,y)}"] & {\c(\lambda a,y)} \arrow[l] \arrow[u, "{\c(\eta_a,y)}", two heads] & {\c(b'',y)} \arrow[l] \arrow[u, "{\c(\kappa_2\kappa_1,y)}"] & {\c(\lambda c,y)} \arrow[l] \arrow[u, "{\c(\eta_c,y)}", two heads] \\
       & {\E(c_{(\calr^{\vee})},y)} \arrow[u]                      &                                                                    &                                                             &                                                                   
\end{tikzcd}
\end{equation*}
The additional exact sequence along the left column is induced by the conflation in \eqref{confyconf}, replacing $x$ with $c$. The map $\E(\eta_c,y)$ is a monomorphism since ${\E(c_{(\calr^{\vee})},y)}$ vanishes. Using the Four Lemma \cite[Lemma I.3.2]{Mac67}, we conclude that $\c(\kappa_2\kappa_1,y)$ is an epimorphism, and in turn that $(\lambda,\psi)$ is an extriangulated functor. More importantly, we deduce that the functor $F\colon \c\to\overline{\c}_{\calr}$ is extriangulated.
	
	Having constructed an {extriangulated} functor $F\colon \c \to \overline{\c}_{\calr}$, we move on to show that it induces a functor $\overline{F}\colon \c/\thick(\calr) \to \overline{\c}_{\calr}$. Given the definition of Verdier quotients in \Cref{NOSloc}, the following is sufficient for this to occur: if $f$ is an inflation with cone in $\thick(\calr)$, or a deflation with cocone in $\thick(\calr)$, then the functor $F\colon \c \to \overline{\c}_{\calr}$ sends $f$ to an isomorphism. Since $F$ is extriangulated, one treats either case by showing that any object in $\thick(\calr)$ is sent to a zero object.	Looking back at \eqref{eq:diadia}, we show that $z\in \calr$ as soon as $x\in \thick(\calr)$. Since $\calr^{\wedge}\subseteq \thick(\calr)$, it follows from extension-closure and \Cref{cotor_facts} that $x^{({^{\perp}\calr})}\in \thick(\calr)\cap {^{\perp}\calr} = \calr^{\vee}$, and in turn that $z\in \calr^{\vee}$, which shows that $\overline{\lambda\rho}(x)\in \calr^{\vee}/[\calr]$. An orthogonal approach shows that $\overline{\rho\lambda}(x)\in \calr^{\wedge}/[\calr]$. Since the functors $\overline{\lambda\rho}$ and $\overline{\rho\lambda}$ are naturally isomorphic to $F$, we conclude that $F(x)\in (\calr^{\vee}\cap \calr^{\wedge})/[\calr] = \calr/[\calr]$, rendering it a zero object. We have thus induced a functor $\overline{F} \colon \c/\thick(\calr) \to \overline{\c}_{\calr}$.
	
	Lastly, we show that $\overline{F} \colon \c/\thick(\calr) \to \overline{\c}_{\calr}$ is a quasi-inverse of the functor $\Phi\colon \overline{\c}_{\calr} \to \c/\thick(\calr)$ induced by the Verdier localization functor $L_{\calr}$. This is a direct consequence of universal properties. Indeed, the functor $\Phi$ the unique functor up to natural isomorphism with the property that $\Phi F\vert_{\z_{\calr}} = L_{\calr}\vert_{\z_{\calr}}$, and $\overline{F}$ is the unique functor up to natural isomorphism with the property that $\overline{F}L_{\calr} = F$, whence $\Phi$ and $\overline{F}$ must be mutually quasi-inverse. Having factored the Verdier localiztion $L_{\calr}$ up to natural isomorphism as $\c \to^{F} \overline{\c}_{\calr} \to^{\Phi} \c/\thick(\calr)$, where $F$ is extriangulated and $\Phi$ is an equivalence, we hereby conclude the proof.
	\end{proof}

\subsection{Example}\label{monica} Let $\Lambda$ be a finite-dimensional $\k$-algebra and let $r$ be a rigid object in the two-term category $\c=H_0\twoper(\Lambda)$. {Then $\c/\thick(\calr)\simeq \overline{\c}_{\add(r)} \simeq H_0\twoper(C_{r})$, where $C_{r}$ is a connective dg $\k$-algebra \cite{Gar24}. 
The finite-dimensional $\k$-algebra $H_0C_{r}$ is the $\tau$-tilting reduction of $\Lambda$ with respect to $r$, in the sense of Jasso \cite{Jas15}.

\subsection{} Although \nref{cor:siltred0Aus}\ref{prop:siltred} shows that the Verdier localization functor $L_{\calr}\colon \c\to \c/\thick(\calr)$ is an additive functor, it is not true in general that the thick subcategory ${\thick(\calr)}$ admits a calculus of left fractions.\footnote{The reader may check that the following is a counter-example: take $\c$ to be the two-term category of the preprojective algebra of type $A_2$, and $\calr$ the additive closure of an indecomposable projective object.}
However, every morphism in $\c/\thick(\calr)$ takes the form $L_{\calr}(s)^{-1} L_{\calr}(f) L_{\calr}(t)^{-1}$, where $s,t\in \mathcal{S}_{\thick(\calr)}$
\cite[Corollary 3.11]{Oga22}, where $\mathcal{S}_{\thick(\calr)}$ is as defined in \Cref{NOSloc}. Since the expression $L_{\calr}(s)^{-1} L_{\calr}(f) L_{\calr}(t)^{-1}$ can be thought of as the gluing of a left fraction and a right fraction, this is reasonably close to a calculus of fractions. 
	 
 \subsection{Corollary}\label{cor:siltred_0Aus_dg}
Let $\A$ be a 0-Auslander exact dg $\k$-category with $\k$-cofibrant mapping complexes.
Fix a rigid dg subcategory $\r\subseteq\A$, and let $\mathscr{Z}_{\r}$ be the dg subcategory of $\A$ such that $H_0\mathscr{Z}_{\r} = \z_{H_0\r}$. Suppose that the rigid subcategory $H_0\r\subseteq H_0\A$ admits Bongartz completions. 
\begin{enumerate}
\item\label{cor:siltred_0Aus_dg0'} The dg quotient $\mathscr{Z}_{\r}/\r$ admits an exact structure such that $H_0(\mathscr{Z}_{\r}/\r)$ is extriangulated equivalent to the ideal quotient $H_0\mathscr{Z}_{\r}/[H_0\r]$.
\item\label{cor:siltred_0Aus_dg0} The inclusion of $\mathscr{Z}_{\r}$ into $\A$ induces a quasi-equivalence $\begin{tikzcd} \mathscr{Z}_{\r}/ {\r} \arrow[r,"\sim"] & \A/\r.  \end{tikzcd}$
\item\label{cor:siltred_0Aus_dg0.5} Let $\zeroAus_\k$ denote the full subcategory of $\Hqe_\k$ spanned by the $0$-Auslander exact dg $k$-categories. The dg quotient $\A/\r$ admits a 0-Auslander exact structure in such a way that the localization morphism $\varpi_{\r}\colon\A\to\A/\r$ becomes an exact morphism in $\zeroAus_{\k}$. 

\item\label{cor:siltred_0Aus_dg1} The localization morphism $\varpi_{\r}\colon\A\to\A/\r$ in $\Hqe_{\k}$ induces a bijection
\begin{equation}\label{eq:rigidred_0Aus_dg}
\varpi_{\r}\colon \rigid_{\r}(\A) \to \rigid(\A/\r),
\end{equation}
that restricts to an isomorphism of posets
\begin{equation}\label{eq:siltingred_0Aus_dg}
\varpi_{\r}\colon \silt_{\r}(\A) \to \silt(\A/\r).
\end{equation}
\item\label{cor:BorComp} Let $\q$ be a rigid dg subcategory of $\A$ such that $\r\subseteq \q$.  Suppose that \nameref{ass:cotor} holds for and $H_0\q$ as a rigid subcategory of $H_0\A$. Then the morphisms $\varpi_{\varpi_{\r}(\q)}\circ \varpi_{\r}$ and $\varpi_{\q}$ coincide in $\Hqe_{\k}$, and thus induce the same map $\rigid_{\r}(\A)\to \rigid(\A/\q)$ (as well as the restriction $\silt_{\r}(\A)\to \silt(\A/\q)$).
\end{enumerate}
\begin{proof}
By \nref{lem:red_0Aus}\ref{lem:red_0Aus2}, we have that $\mathscr{Z}_{\r}$ is a 0-Auslander exact dg $\k$-category. Moreover, it follows directly from the definitions that all objects in $H_0\r$ are projective-injective in $H_0\mathscr{Z}_{\r}$. The statement in \ref{cor:siltred_0Aus_dg0'} now follows from work of Chen \cite[Theorem 6.35 and its proof]{Che23thesis}.

We will now prove \ref{cor:siltred_0Aus_dg0}. The inclusion of $\mathscr{Z}_{\r}$ into $\A$ clearly induces quasi-isomorphisms of complexes $\mathscr{Z}_{\r}(x,y)\simeq\A(x,y)$ for any pair of objects $x,y\in \mathscr{Z}_{\r}$. Passing to the Drinfeld dg quotients, we also have quasi-isomorphisms ${\mathscr{Z}_{\r}\over {\r}}(x,y)\simeq{\A\over {\r}}(x,y)$, which verifies that one of the two properties of a quasi-equivalence holds. The other property to verify is that $\begin{tikzcd} H_0(\mathscr{Z}_{\r}/ {\r})\arrow[r,"\sim"] & H_0(\A/\r) \end{tikzcd}$ is an equivalence of categories. We deduce from \nref{cor:siltred0Aus}\ref{prop:siltred} the Verdier localization functor $L_{\calr}$ is additive, whence \nref{prop:dgquot_ex} yields an equivalence $H_0(\A/\r)\simeq H_0(\A)/\thick(H_0\r)$, and we may then conclude thanks to \nref{cor:siltred0Aus}\ref{prop:siltred} and \ref{cor:siltred_0Aus_dg0'} above. We now use \nref{lem:red_0Aus}\ref{lem:red_0Aus2} to deduce \ref{cor:siltred_0Aus_dg0.5}, and \nref{cor:siltred0Aus}\ref{cor:siltred0Aus1} to deduce \ref{cor:siltred_0Aus_dg1}.

Lastly, we prove \ref{cor:BorComp}.
Using the universal property of dg localizations in \Cref{dgloc}, one shows that the localization morphism $\varpi_{\q}\colon \A\to \A/\q$ in $\Hqe_{\k}$ admits a factorization
	\begin{equation*}
	\begin{tikzcd}
		\A \arrow[r,"\varpi_{\r}"] & \A/\r \arrow[r,"\varpi_{\varpi_{\r}(\q)}"] & (\A/\r)/\varpi_{\r}(\q) \arrow[r,"\simeq"] &\A / \q.
	\end{tikzcd} 
	\end{equation*}
	It follows that the maps $\varpi_{\varpi_{\r}(\q)}\circ \varpi_{\r}$ and $\varpi_{\q}$ coincide. By \ref{cor:siltred_0Aus_dg1} above, the morphisms $\varpi_{\varpi_{\r}(\q)}\circ \varpi_{\r}$ and $\varpi_{\q}$ induce the same map $\rigid_{\r}(\A)\to \rigid(\A/\q)$.
\end{proof}

\section{Picture categories and picture groups}\label{subsec:tcmc}

 \subsection{Definition}\label{subsec:BorComp2}
Let $\A$ be a 0-Auslander exact dg $\k$-category with $\k$-cofibrant mapping complexes. Assume also that all rigid subcategories of $H_0\A$ admit Bongartz completions. If these assumption hold, we say that the \textit{picture category $\tcmc{\A}$ of $\A$ is defined}.
The \textit{picture category} of $\A$ is denoted $\tcmc{\A}$, and is defined in the paragraph below. By \nref{cor:siltred_0Aus_dg}\ref{cor:siltred_0Aus_dg0.5} and \nref{lem:red_0Aus}\ref{lem:red_0Aus3}, the assumptions we have made on $\A$ also hold for dg quotient $\A/\r$ where $\r$ is a rigid dg subcategory of $\A$.

Two dg quotients $\L_1$ and $\L_2$ of $\A$ are \textit{equivalent} if the localization morphisms $\A\to^{\varpi_1} \L_1$ and $\A\to^{\varpi_2} \L_2$ are isomorphic as objects in the under-category (aka. co-slice category) of $\A$ in $\Hqe_{\k}$.
The objects in $\tcmc{\A}$ will be taken as the equivalence classes of dg quotients of $\A$ with respect to rigid dg subcategories. We will often treat an object in $\tcmc{\A}$ as a dg quotient, by implicitly choosing a representative of the equivalence class. By \nref{cor:siltred_0Aus_dg}\ref{cor:siltred_0Aus_dg0.5}, such dg quotients are in $\zeroAus_{\k}$.
A morphism from $\L_1$ to $\L_2$ in $\tcmc{\A}$ is given by a rigid dg subcategory $\r$ of $\L_1$ such that $\L_1/\r \simeq \L_2$ in $\zeroAus_{\k}$.
The morphism in $\tcmc{\A}$ given by $\r$ will be denoted $\L_1\to^{\r}\L_2$. 
We define the composite of $\L_1\to^{{\r_1}}\L_2 \to^{\r_2}\L_3$ by 
$\begin{tikzcd}[column sep=3em]
\L_1\arrow[r,"\varpi^{-1}_{\r_1}(\r_2)"] & \L_3,
\end{tikzcd}$
 where $\varpi_{\r_1}$ is the map defined in \nref{cor:siltred_0Aus_dg}\ref{cor:siltred_0Aus_dg1}. 
 
 \subsection{}\label{pre:tcmc_isCat}
We will now show that picture categories indeed are categories. The first step is to reassure the reader that the composition rule is well-defined. In precise terms, we should confirm that the composite $\begin{tikzcd}[column sep=3em]
\L_1\arrow[r,"\varpi^{-1}_{\r_1}(\r_2)"] & \L_3,
\end{tikzcd}$ of $\L_1\to^{{\r_1}}\L_2 \to^{\r_2}\L_3$ satisfies the following:
\begin{enumerate}
	\item\label{pre:tcmc_isCat1} the subcategory $\varpi^{-1}_{\r_1}(\r_2)$ of $\L_1$ is rigid,
	\item\label{pre:tcmc_isCat2} the dg quotient $(\L_1/\r_1)/\r_2$ is naturally identified with $\L_1/ \varpi_{\r_1}^{-1}(\r_2)$.
\end{enumerate}
The first assertion follows directly from \nref{cor:siltred_0Aus_dg}\ref{cor:siltred_0Aus_dg1}. The second is deduced from \nref{cor:siltred_0Aus_dg}\ref{cor:BorComp}; the dg $\k$-category $\L_1/\r_1$ has $\k$-cofibrant mapping complexes, so the morphisms $\varpi_{\varpi_{\r_1}^{-1}(\r_2)}$ and $\varpi_{\r_2}\circ\varpi_{\r_1}$ coincide in $\zeroAus_{\k}$, which amounts to that $(\L_1/\r_1)/\r_2 \simeq \L_1/ (\varpi_{\r_1}^{-1}(\r_2))$ in $\zeroAus_{\k}$. 

 \subsection{Proposition}\label{prop:tcmc_isCat}
	Let $\A$ be a 0-Auslander exact dg $\k$-category for which the picture category $\tcmc{\A}$ of $\A$ is defined.
	Then $\tcmc{\A}$ is a category.\footnote{The proof below closely resembles that of the special case of connective dg $\k$-algebras with finite-dimensional homology in all degrees, as was established in previous work \cite[Theorem 4.3]{Bor21}}
\begin{proof}
It is clear from the construction that the identity morphism at an object $\L\in\tcmc{\A}$ is provided by the rigid dg subcategory of $\L$ consisting of the zero objects. It remains to show that the composition rule defined above is associative. Given three composable morphisms 
\begin{equation*}
\begin{tikzcd}
	\L_1\arrow[r,"\r_1"] & \L_2 \arrow[r,"\r_2"] & \L_3 \arrow[r,"\r_3"] & \L_4,
\end{tikzcd}
\end{equation*}
it is to be shown that 
\begin{equation}\label{eq:tcmc_isAss}
	\varpi^{-1}_{\varpi_{\r_1}^{-1}(\r_2)}(\r_3) = \varpi^{-1}_{\r_1}\big(\varpi_{\r_2}^{-1}(\r_3)\big),
\end{equation}
as rigid dg subcategories of $\A$. Since $\varpi_{\r_1}^{-1}(\r_2)$ is a rigid dg subcategory of $\A$ containing $\r_1$, \nref{cor:siltred_0Aus_dg}\ref{cor:BorComp} enables the decomposition of the localization morphism $\varpi_{\varpi_{\r_1}^{-1}(\r_2)}$ as $\varpi_{\r_2}\circ \varpi_{\r_1}$. As a result, we have that $\varpi_{\varpi_{\r_1}^{-1}(\r_2)}$ and $\varpi_{\r_2}\circ \varpi_{\r_1}$ induce the same bijection from $\rigid_{\varpi_{\r_1}^{-1}(\r_2)}(\L_1)$ to $\rigid(\L_3)$.
The application of the inverses of these coincident bijections to $\r_3$ gives the equality in \eqref{eq:tcmc_isAss}, whence the proof is complete.
\end{proof}

 \subsection{Definition}\label{def:pic}
	Let $\A$ be a 0-Auslander exact dg $\k$-category for which the picture category $\tcmc{\A}$ of $\A$ is defined. 
		\begin{enumerate}
		\item\label{def:pic1} The \textit{picture space} of $\A$ is denoted $\picspace{\A}$, and defined as the classifying space $\classspace \tcmc{\A}$. 
		\item\label{def:pic2} The \textit{picture group} of $\A$ is denoted $\picgroup{\A}$, and defined as the fundamental group $\pi_1 (\picspace{\A},\mathscr{O})$, where $\mathscr{O}$ is a trivial dg $\k$-category. 
	\end{enumerate}
	
\subsection{Example}\label{eg:Bor21}
	Let $A$ be a dg $k$-algebra with finite-dimensional homology in all degrees and let $\A=\twoper(A)$ be the two-term dg $\k$-category of $A$  (see \nref{eg:0Aus}\ref{eg:0Aus_per2}). Then the category $\tcmc{\A}$ is equivalent to the $\tau$-cluster morphism category of $A$, as defined in previous work \cite{Bor21}. This can be shown directly by spelling out the definition provided in the cited reference. 
	For this class of examples, we simplify the notation, writing $\tcmc{A}$ for $\tcmc{\twoper(A)}$, and so on. 
			
The dg $\k$-algebra homomorphism $A\to H_0A$ induces an equivalence of picture categories $\tcmc{\A} \to^{\sim} \tcmc{H_0A}$ \cite[Theorem 4.4]{Bor21}, whence one can restrict to the special case of finite-dimensional $\k$-algebras without loss of generality, up to equivalence of categories.
For finite-dimensional $\k$-algebras, there are several instances of worked examples in the literature \cite[§12]{BM18w} \cite[§7]{BH21}. A particularly easy example is that of a local $\k$-algebra $\Lambda$, whose picture category can be displayed as 
\begin{tikzcd}
\twoper(\Lambda) \arrow[r,yshift=0.25em,"\Lambda"] \arrow[r,yshift=-0.25em,"\Sigma\Lambda"'] & \mathscr{O},
\end{tikzcd}
where the arrows are labeled by the rigid object generating the rigid dg subcategory of $\twoper(\Lambda) $ in question. The picture space of $\Lambda$ is a circle, and the picture group becomes free abelian with one generator.

We include a small example to illustrate what can happen in the non-reduced case. The 0-Auslander extriangulated $\k$-category $\mathrm{mod}(\Lambda_2)$ in \nref{eg:0Aus}\ref{eg:0Aus_mod} is Quillen exact and contains a projective generator, whence it admits an exact dg enhancement $\A$ whose picture category is defined. The $\k$-category $\mathrm{mod}(\Lambda_2)$ contains three objects up to isomorphism, namely $n$, $p$, and $i$, where $n$ is projective-injective, the object $p$ is projective and non-injective, and $i$ is injective and non-projective. 
Every indecomposable object in $\mathrm{mod}(\Lambda_2)$ is rigid, and will be tacitly identified with the rigid dg subcategory of $\A$ it generates. We thus write $\A/x$ for the dg quotient of $\A$ modulo the smallest additive dg subcategory of $\A$ containing $x$.  
The picture category $\tcmc{\A}$ can then be displayed as follows:
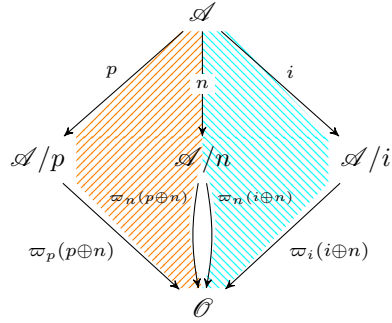
\begin{equation*}
	\begin{tikzcd}[row sep=4em, column sep=1.62em,
	ampersand replacement=\&,
	execute at end picture={
    \scoped[on background layer]
    \fill[pattern=north east lines, pattern color=orange] (a.south west) -- (ap.north east) -- (ap.south east) -- (o.north west) -- (o.north) to [ bend left=20]  (an.center) -- (a.south) -- cycle;},
    execute at end picture={
    \scoped[on background layer]
    \fill[pattern=north west lines, pattern color=cyan] (a.south east) -- (ai.north west) -- (ai.south west) --  (o.north east) -- (o.north) to [ bend right=20]  (an.center) -- (a.south) -- cycle;}]
                                               \&        \&  |[alias=a]| \A \arrow[lld, "p"'] \arrow[d, "n" description] \arrow[rrd, "i"]                           \&        \&                                               \\
 |[alias=ap]|\A/p \arrow[rrd, "{\varpi_{p}(p\oplus n)}"'] \&  \&  |[alias=an]| \A/n \arrow[d, "{\scriptscriptstyle\varpi_{n}(p\oplus n)}"', near start, bend right=10] \arrow[d, "{\scriptscriptstyle\varpi_{n}(i\oplus n)}", near start, bend left=10] \& \&  |[alias=ai]| \A/i \arrow[lld, "\varpi_{i}(i\oplus n)"] \\
                                               \&        \&  |[alias=o]|\mathscr{O}                                                                                  \&        \&                                              
\end{tikzcd}
\end{equation*}
where the two shaded squares are commutative. There are two morphisms from $\A$ to $\mathscr{O}$, namely $p\oplus n$ and $i\oplus n$.
Note that $\A/n$ is quasi-equivalent to $\twoper(\k)$, and that the full subcategory of $\tcmc{\A}$ spanned by $\A/n$ and $\mathscr{O}$ is equivalent to $\tcmc{\k}$. It is moreover easy to see that the picture space $\picspace{\A}$ to homotopy equivalent to $\picspace{\k}$, an assertion we generalize in \nref{prop:homoeq} below.

\subsection{Proposition}\label{prop:homoeq}
{Let $\A$ be a 0-Auslander exact dg $\k$-category for which the picture category $\tcmc{\A}$ of $\A$ is defined, and let $\mathscr{N}$ denote the rigid dg subcategory of $\A$ consisting of all projective-injective objects in $\A$. The space} $\picspace{\A/\mathscr{N}}$ {is then a deformation retract of $\picspace{\A}$.
In particular, the inclusion functor $\begin{tikzcd}\tcmc{\A/\mathscr{N}} \arrow[r,hook] & \tcmc{\A}\end{tikzcd}$ induces a homotopy equivalence $\begin{tikzcd}\picspace{\A/\mathscr{N}} \arrow[r,"\sim"] & \picspace{\A}\end{tikzcd}$ and thus an isomorphism of picture groups $\begin{tikzcd}\picgroup{\A/\mathscr{N}} \arrow[r,"\sim"] & \picgroup{\A}.\end{tikzcd}$
{One consequently recovers all possible homotopy types of picture spaces, as well as all isomorphism classes of picture groups, by only considering reduced 0-Auslander exact dg $\k$-categories.}
More generally, if $\r$ is a rigid dg subcategory of $\A$ in which all objects are projective-injective, the inclusion functor $\begin{tikzcd}\tcmc{\A/\r} \arrow[r,hook] & \tcmc{\A}\end{tikzcd}$  induces a homotopy equivalence of picture spaces} $\begin{tikzcd}\picspace{\A/\r} \arrow[r,"\sim"] & \picspace{\A}\end{tikzcd}$ and thus an isomorphism of picture groups $\begin{tikzcd}\picgroup{\A/\r} \arrow[r,"\sim"] & \picgroup{\A}.\end{tikzcd}$
\begin{proof}
Since $\r$ is a full dg subcategory of $\mathscr{N}$, we have consecutive inclusions 
\begin{equation*}
\begin{tikzcd} \picspace{\A/\mathscr{N}} \arrow[r,hook] & \picspace{\A/\r} \arrow[r,hook] & \picspace{\A}. \end{tikzcd}
\end{equation*}
One uses the two-out-of-three property for homotopy equivalences to deduce our claim concerning $\r$ from our claim concerning $\mathscr{N}$. To complete the proof, we show that the inclusion functor $\begin{tikzcd}\tcmc{\A/\mathscr{N}} \arrow[r,hook] & \tcmc{\A}\end{tikzcd}$ admits a left adjoint \cite[Corollary 1 in §1]{Qui73}. 
For an object $\L\in \tcmc{\A}$, let $\mathscr{N}_{\L}\subseteq \L$ denote the rigid dg subcategory consisting of the projective-injective objects. We will define an endofunctor $\lambda$ on $\tcmc{\A}$ sending $\L$ to its reduced counterpart, i.e. to $\L/\mathscr{N}_{\L}$ (cf. \Cref{redcat}). If $\L \to^{\r} \L'$ is an arbitrary morphism in $\tcmc{\A}$, we can, as a result of \nref{prop:dgquot_ex}\ref{prop:dgquot_ex2} and \nref{cor:siltred0Aus}\ref{prop:siltred}, regard $H_0\L'$ as an ideal quotient of a subcategory $\z_{H_0\r}$ of $H_0\L$ containing all projective-injectives. Consequently, the extriangulated functor $H_0\L \to H_0\L'$ sends projective-injectives to projective-injectives, whence we may induce a morphism $\lambda\L \to^{} \lambda\L'$ in $\tcmc{\A}$. It is now easy to see that $\lambda$ is a functor. Moreover, the image of $\lambda$ is precisely $\tcmc{\A/\mathscr{N}}$, since any morphism of the form $\A \to \lambda\L$ in $\tcmc{\A}$ factors through $\A\to^{{\mathscr{N}}} \A/\mathscr{N}$. We will conclude by showing that $\begin{tikzcd}\lambda\colon\tcmc{\A} \arrow[r] & \tcmc{\A/\mathscr{N}}\end{tikzcd}$is left adjoint to the inclusion functor $\begin{tikzcd}\tcmc{\A/\mathscr{N}} \arrow[r,hook] & \tcmc{\A};\end{tikzcd}$ indeed, given a morphism $\L \to \L'$ in $\tcmc{\A}$, where $\L'\in \tcmc{\A/\mathscr{N}}$, it factors uniquely through $\L \to \L/\mathscr{N}_{\L}$, providing a natural bijection
\begin{equation*}
	\tcmc{\A}(\L,\L') \cong \tcmc{\A/\mathscr{N}}(\lambda\L,\L').
\end{equation*}
as desired.
\end{proof}

\subsection{}\label{tcmc_str} 
Henceforth in this section, we describe the structure of picture categories and picture groups. Note first that $\tcmc{\A}$ is a skeletal category by definition, in the sense that isomorphic objects in $\tcmc{\A}$ are equal. 
The object $\A$ in $\tcmc{\A}$ is \textit{weakly initial}, i.e. there exists a morphism from $\A$ to any other object in $\tcmc{\A}$ by definition. The trivial dg $\k$-category $\mathscr{O}$ is a {weakly terminal object} in $\tcmc{\A}$, since it occurs as a dg quotient of any $\L\in\tcmc{\A}$ modulo a silting subcategory of $\L$. The existence of a weakly initial (or weakly terminal) object in $\tcmc{\A}$  makes picture space $\picspace{\A}$ a connected topological space.
The morphism set $\tcmc{\A}(\A,\mathscr{O})$ is precisely the set of silting subcategories of $\A$. More generally, the morphism set $\tcmc{\A}(\L,\mathscr{O})$ is precisely the set of silting subcategories of $\L$. Unless $\L$ is equal to $\A$ (resp. $\mathscr{O}$), the morphism set $\tcmc{\A}(\L,\A)$ (resp. $\tcmc{\A}(\mathscr{O},\L)$) will be empty. 
For a rigid dg subcategory $\r$ of $\A$, the picture category $\tcmc{\A/\r}$ embeds into $\tcmc{\A}$. 

\subsection{}\label{prop:HI18.2.14}
 Let $\A$ be a 0-Auslander dg $\k$-category for which the picture category $\tcmc{\A}$ of $\A$ is defined. Suppose that $H_0\A$ is Krull--Schmidt and contains a silting object $t$ with $n$ non-isomorphic indecomposable direct summands.  By combining the assertions in the paragraphs below, we will show that $\tcmc{\A}$ is a \textit{cubical category}, in the sense of Igusa \cite[Definition 3.2]{Igu14}, whereby generalizing a result of Hanson--Igusa, who show that the picture category of a finite-dimensional $\k$-algebra is cubical \cite[Theorem 2.14]{HI21}.
 
 Let $\L_1 \to^{\r} \L_2$ be a morphism in $\tcmc{\A}$. By \Cref{AT.5.4}, 
 the rigid subcategory $H_0\r$ of $H_0\L_1$ is of the form $\add(r)$ for some rigid object $r$ of $H_0\L_1$, which can be assumed to be basic (i.e. its indecomposable objects are pair-wise non-isomorphic). We define the \textit{rank} of $\r$ as the number of non-isomorphic indecomposable direct summands of $r$. By \Cref{AT.5.4} again, the rank of a morphism is bounded by $n$, whence it is a well-defined number. Using the definition of the composition rule in $\tcmc{\A}$ directly, one sees that the rank of a composite morphism $\L_1 \to^{\r_1} \L_2 \to^{\r_2}\L_3$ is the sum of the ranks of $\r_1$ and $\r_2$.
 
 Let $\L_1 \to^{\r} \L_2$ be a morphism in $\tcmc{\A}$ of rank $\ell$. Let $r$ be a basic rigid object of $H_0\L_1$ such that $H_0\r = \add(r)$ as full subcategories of $H_0\L_1$.
 The rigid object $r\in H_0\L_1$ can then be written as a biproduct $r \simeq \bigoplus_{i=1}^\ell r_i$, where indecomposable objects $r_i$ are pair-wise non-isomorphic. If we factor $\r$ as $\L_1 \to^{\r'} \L_{1.5} \to{} \L_2$, then $\r'$ must be given by a direct summand of $r$. Let $\mathrm{faq}(\r)$ denote the subcategory of $\tcmc{\A}$ consisting of the morphisms $\mathscr{Q}$ such that $\r = \mathscr{Q}_2 \circ \mathscr{Q} \circ \mathscr{Q}_1$ for some morphisms $\mathscr{Q}_1$ and $\mathscr{Q}_2$. One shows that $\mathrm{faq}(\r)$ is the image of a faithful embedding into $\tcmc{\A}$ of a Boolean lattice of rank $\ell$ (there is one object in $\mathrm{faq}(\r)$ for every direct summand of $r$). This is to say that $\mathrm{faq}(\r)$ takes the form of an $\ell$-dimensional hypercube.
 
We keep the notation introduced in the last paragraph.
The set $\{\r_i\}$ of rigid dg subcategories of $\L_1$ provides a set of morphisms $\{\L_1 \to^{\r_i} \L_i \}$ (where the objects $\L_i$ need not be pair-wise distinct). The elements of this set are called the \textit{first factors} of $\r$. A morphism is uniquely determined by its set of first factors, simply because the rigid object $r$ is uniquely determined by its indecomposable direct summands. The set of \textit{last factors} of $\r$ is given by {$\{\L_j \to^{\varpi_{\r^i}(\r)} \L_2 \}$, where $\r^i$ is the dg subcategory such that $H_0(\r^i) = \add(r/r_i)$}, and this set also uniquely determines the morphism $\r$.

\subsection{} 
Although we have defined picture groups as the fundamental group of the picture space, we will in practice treat it as the fundamental group of the picture category. We slightly adapt the setup of Evrard \cite{Evr73,Evr75} \cite[§1]{Qui73}. A \textit{loop} in $\tcmc{\A}$ with basepoint $\mathscr{O}$ is either the trivial path on $\mathscr{O}$ or a zig-zag $\gamma$ of morphisms in $\tcmc{\A}$ of the form
	\begin{equation}\label{eq:EvrLoop}
		\begin{tikzcd}[row sep=0.6em]
 \mathscr{O} & \arrow[l] \L_m \arrow[r] & {\L_{m-1,m}} & \cdots \arrow[l] \arrow[r] & \L_2 \arrow[r] & {\L_{1,2}} & \L_1 \arrow[r] \arrow[l] & \mathscr{O}
\end{tikzcd}
	\end{equation}
	 The \textit{length} of the loop above is $m$, whereas the trivial loop has length 0. In \eqref{eq:EvrLoop}, the object $\L_i$ is the \textit{$i^{\mathit{th}}$ term} of $\gamma$.
	 
	Let $\gamma_1$ and $\gamma_2$ be loops on $\mathscr{O}$, of length $m$ and $p$, respectively. 
	The \textit{concatenation} of $\gamma_1$ and $\gamma_2$ is denoted by $\gamma_2\gamma_1$ and defined by 
	\begin{equation*}
		\begin{tikzcd}[row sep=0.5em]
                                      &                            &        &                             & \mathscr{O}                           &                             &        &                            &             \\
\mathscr{O}                           & \L_p^2 \arrow[l] \arrow[r] & \cdots & \L_1^2 \arrow[l] \arrow[ru] &                                       & \L_m^1 \arrow[lu] \arrow[r] & \cdots & \L_1^1 \arrow[l] \arrow[r] & \mathscr{O} \\
{} \arrow[rrrr, ""',symbol=\underbrace{\hspace{17.5em}}_{\gamma_2}] &                            &        &                             & {} \arrow[rrrr, ""', symbol=\underbrace{\hspace{17.5em}}_{\gamma_1}] &                             &        &                            & {}       
\end{tikzcd}	\end{equation*}
	An \textit{elementary homotopy} from  $\gamma_1$ to $\gamma_2$ is given by a strictly increasing function $f\colon [1,m]\cap \Z \rightarrow [1,p]\cap \Z$ and a commutative diagram
	\begin{equation}\label{eq:elhom}
	\begin{tikzcd}
f(\gamma_1)\colon & \mathscr{O} \arrow[d, equal] & \L_m^1 \arrow[l] \arrow[r] \arrow[d] & \cdots \arrow[r] & {\L_{1,2}^1} \arrow[d] & \L_{1}^1 \arrow[l] \arrow[d] \arrow[r, equal] & \cdots & \L_{1}^1 \arrow[l, equal] \arrow[d] \arrow[r] & \mathscr{O} \arrow[d, equal] \\
\gamma_2\colon    & \mathscr{O}                    & \L_p^2 \arrow[r] \arrow[l]                    & \cdots \arrow[r] & \L_{f(1),f(1)+1}^2                     & \L_{f(1)}^2 \arrow[l] \arrow[r]                        & \cdots & \L_{1}^2 \arrow[l] \arrow[r]                             & \mathscr{O}                   
\end{tikzcd}
	\end{equation}
	where the bottom row displays $\gamma_2$ and the top row the \textit{elongation} of $\gamma_1$ along $f$ (there are $f_{i+1}-f_{i}$ copies of $\L_i$ in the top row, where $f_i\defeq f(i)$ for $i\geq 1$ and $f_0\defeq 0$). 
	The notion of elementary homotopy determines a binary relation on the set of loops, whose symmetric and transitive closure provides a suitable notion of homotopy. Explicitly, two loops are said to be \textit{homotopic} if there exists a zig-zag of elementary homotopies from one to the other.
	
	With the definitions above, one defines the fundamental group of $\tcmc{\A}$ as the set homotopy classes of loops on $\mathscr{O}$ equipped with the binary operation of concatenation. This is group is clearly isomorphic to the picture group of $\A$, as defined in \nref{def:pic}\ref{def:pic2}.
	
	\subsection{} Recall that an \textit{interval} in a poset $({P},\leq)$ is a subposet of ${P}$ of the form $$[a,b]\defeq\{x\in {P} \sth a\leq x\leq b \}, $$ where $a,b\in {P}$. Let $\mathrm{int}({P})$ denote the set of intervals in ${P}$. An interval $[a,b]$ is said to be \textit{irreducible} if its cardinality is exactly $2$.
	
Our next result should be compared with presentations of picture groups that have been established for the special case of $\bfg$-finite finite-dimensional algebras over fields \cite[Proposition 4.4]{HI21}. 

\subsection{Proposition}\label{prop:picgroup_pres}
Let $\A$ be a 0-Auslander exact dg $\k$-category for which the picture category $\tcmc{\A}$ of $\A$ is defined. Let $\silt(\A)$ be the poset of silting subcategories of $\A$, as defined in \Cref{AT22.5.12}. Consider the (not necessarily injective) map $\gamma_{(-)} \colon \mathrm{int}(\silt(\A)) \to \picgroup{\A}$ mapping an interval $[\s_1, \s_2]\subseteq \silt(\A)$ to the homotopy class of the loop
\begin{equation}\label{eq:picgroup_pres_genloop}
\begin{tikzcd}[column sep=3em]
			 \mathscr{O} & \arrow[l,"\s_2"'] \A \arrow[r,"\s_1"] & \mathscr{O}. 
		\end{tikzcd}
\end{equation}
\begin{enumerate}
	\item\label{prop:picgroup_pres1} The image of $\gamma_{(-)}$ generates $\picgroup{\A}$ as a group. In particular, the map $\gamma_{(-)}$ induces a surjective homomorphism of groups
	\begin{equation}\label{eq:picgroup_pres1}
		{\gamma_{(-)}}\colon  {G_{\mathrm{int}(\silt(\A))}} \epi \picgroup{\A},
	\end{equation}
	where ${G_{\mathrm{int}(\silt(\A))}}$ is the free group generated by ${\mathrm{int}(\silt(\A))}$. Thus, if $\A$ is $\bfg$-finite, the picture group $\picgroup{\A}$ is finitely generated. 
	\item\label{prop:picgroup_pres2} The kernel of the homomorphism \eqref{eq:picgroup_pres1} is generated by expressions of the form  $${[\s_1,\s_2]}^{-1}{[\s_1, \s_{1.5}]}{[\s_{1.5}, \s_{2}]},$$ where $[\s_1, \s_{1.5}]$ and $[\s_{1.5}, \s_2]$ are consecutive intervals in $\silt(\A)$. In particular, if $\A$ is $\bfg$-finite, the picture group $\picgroup{\A}$ is finitely presented, and generated by the irreducible intervals in $\silt(\A)$.
	\item\label{prop:picgroup_pres2.125} Suppose that $\A$ is $\bfg$-finite, and that $H_0\A$ is a Krull--Schmidt extriangulated $\k$-category containing a silting object.
	Let $\mathfrak{B}(\A)$ denote the set of objects $\L\in \tcmc{\A}$ such that every non-zero rigid subcategory of $\L$ is silting and $|\tcmc{\A}(\L,\mathscr{O})|=2$, and let $G_{\mathfrak{B}(\A)}$ be the free group generated by ${\mathfrak{B}(\A)}$.
 	Writing $\tcmc{\A}(\L,\mathscr{O}) = \{\p_{\L}, \i_{\L} \}$, where $\p_\L\geq \i_{\L}$ in $\silt(\L)$, we have that the group homomorphism 
	\begin{equation*}
		\overline{\gamma}_{(-)}\colon G_{\mathfrak{B}(\A)} \to \picgroup{\A},
	\end{equation*}
	sending $\L\in \mathfrak{B}(\A)$ to the homotopy class of the loop 
	\begin{tikzcd}[column sep=3em]
			 \mathscr{O} & \arrow[l,"\p_{\L}"'] \L \arrow[r,"\i_{\L}"] & \mathscr{O},
		\end{tikzcd}
	is surjective.
\end{enumerate}
\begin{proof}
	We address \ref{prop:picgroup_pres1} by showing that an arbitrary loop on $\mathscr{O}$ is homotopic to a concatenation of loops that are either in the image of $\gamma_{(-)}$, or inverses of such. We proceed by induction on the length of loops, treating 0 and 1 as anchor steps. For the loop of length $0$, the claim is trivial is establish. 
	A loop of length 1 is shown along the top row in \eqref{prop:picgroup_pres_loop1} below.
	\begin{equation}\label{prop:picgroup_pres_loop1}
		\begin{tikzcd}
			 \mathscr{O} &  \arrow[l,"\q_1"']  \L_1  \arrow[r,"\q_2"] & \mathscr{O} \\
			 \mathscr{O} \arrow[u,equal] &  \arrow[l,"\q'_1"'] \A \arrow[u,"\r"]  \arrow[r,"\q'_2"] & \mathscr{O}\arrow[u,equal]
		\end{tikzcd}
	\end{equation}
	Since $\A$ is weakly initial in $\tcmc{\A}$, one can pick a morphism $\A \to^{\r} \L_1$ and use it to construct the elementary homotopy shown in \eqref{prop:picgroup_pres_loop1}.
	Consider the maximal element $\mathscr{P}_{\A}$ in $\silt(\A)$ (see \Cref{AT.5.4}). We can decompose the loop along the bottom row in \eqref{prop:picgroup_pres_loop1} as
	\begin{equation*}
		\begin{tikzcd}
			 \mathscr{O} &  \arrow[l,"\q'_1"'] \A  \arrow[r,"\mathscr{P}_{\A}"] & \mathscr{O} & \arrow[l,"\mathscr{P}_{\A}"'] \A \arrow[r,"\q'_2"]  & \mathscr{O}.
		\end{tikzcd}
	\end{equation*}
	Having shown that the loop along the top row of \eqref{prop:picgroup_pres_loop1} is homotopic to $\gamma_{[\q'_1,\mathscr{P}_{\A}]}^{-1} \gamma_{[\q'_2,\mathscr{P}_{\A}]}$, we have settled the anchor steps.
	For the inductive step, suppose that the claim has been shown for all loops of length shorter than $m$, for some positive number $m\geq 2$, and consider an arbitrary loop
		\begin{equation*}
		\begin{tikzcd}
\omega\colon & \mathscr{O} & \L_m \arrow[l] \arrow[r] & \cdots & \L_2 \arrow[l] \arrow[r] & {\L_{1,2}} & \L_1 \arrow[l, "\q_2"'] \arrow[r, "\q_1"] & \mathscr{O}
\end{tikzcd}
	\end{equation*}
	of length $m$.
	Since $\A$ is weakly initial in $\tcmc{\A}$, there exists a morphism $\A \to^{\r_1} \L_1$. Similarly, the weak terminality property of $\mathscr{O}$ lets us choose a morphism $\L_{1,2} \to^{\r_2} \mathscr{O}$ in $\tcmc{\A}$. One sets up a homotopy 
	\begin{equation*}
		\begin{tikzcd}
\mathscr{O}                                       & \L_m \arrow[r] \arrow[l]                                       & \cdots & \L_2 \arrow[r] \arrow[l]                                       & {\L_{1,2}}                                               & \L_1 \arrow[r, "\q_1"] \arrow[l, "\q_2"']                                        & \mathscr{O}                                       \\
\mathscr{O} \arrow[u, equal] \arrow[d, equal] & \L_m \arrow[r] \arrow[l] \arrow[u, equal] \arrow[d, equal] & \cdots & \L_2 \arrow[r] \arrow[l] \arrow[u, equal] \arrow[d, equal] & {\L_{1,2}} \arrow[u, equal] \arrow[d, "\r_2"] & \A \arrow[r, "\s_1"] \arrow[l, "\q_2\circ\r_1"'] \arrow[u, "\r_1"'] \arrow[d, equal] & \mathscr{O} \arrow[u, equal] \arrow[d, equal] \\
\mathscr{O}                                       & \L_m \arrow[r] \arrow[l]                                       & \cdots & \L_2 \arrow[r] \arrow[l]                                       & \mathscr{O}                                              & \A \arrow[r, "\s_1"] \arrow[l, "\s_2"']                                          & \mathscr{O}                                      
\end{tikzcd}
	\end{equation*} 
	where $\s_1=\q_1\circ\r_1$ and $\s_2=\r_2\circ \q_2 \circ \r_1$. The bottom row is now a composite of loops of length shorter than $m$. By the induction hypothesis, we can express the bottom row as a concatenation of loops in the image of $\gamma_{(-)}$, or inverses of such, whereby we conclude the proof of \ref{prop:picgroup_pres1}.
		
We move on to \ref{prop:picgroup_pres2}. It is clear that the kernel of the homomorphism $\gamma$ is generated by expressions of the form $\gamma_2^{-1} \gamma_1$, where  $\gamma_1$ and $\gamma_2$ are related by elementary homotopy. Using \ref{prop:picgroup_pres1}, we can further assume said elementary homotopy to be of the following form
\begin{equation*}\label{prop:picgroup_pres_genloop}
\begin{tikzcd}[column sep=3em]
\mathscr{O} \arrow[d,equal] & \cdots \arrow[l] \arrow[r] & \mathscr{O} \arrow[d] & \A \arrow[d,equal] \arrow[r, equal] \arrow[l, "\s_1"'] & \A \arrow[d, "\s_{1.5}"] & \A \arrow[r, "\s_2"] \arrow[d,equal] \arrow[l, equal] & \mathscr{O} \arrow[d,equal] \\
\mathscr{O}           &  \cdots \arrow[r] \arrow[l]  & \mathscr{O}           & \A \arrow[r, "\s_{1.5}"] \arrow[l, "\s_1"']        & \mathscr{O}              & \A \arrow[r, "\s_2"] \arrow[l, "\s_{1.5} "']      & \mathscr{O}           \\
                     &     {}                      &   \arrow[rrrr,symbol=\underbrace{\hspace{20em}}_{D}]   &   {}   & {}  & {}  &   {}                 
\end{tikzcd}
\end{equation*}
where either $\s_1\leq \s_{1.5} \leq \s_2$ or $\s_2\leq \s_{1.5} \leq \s_1$. The whole diagram can be built from subdiagrams of the same form as $D$, and $D$ is itself an elementary homotopy of loops. Since $D$ encodes a relation of the appropriate form, or the inverse of such an expression, we can conclude that \ref{prop:picgroup_pres2} holds.

We now prove \ref{prop:picgroup_pres2.125}. Having shown in \ref{prop:picgroup_pres2} that $\gamma$ sends the irreducible intervals in $\silt(\A)$ to a set of generators, it suffices to show that every such generator is homotopic to a loop of the form $\overline{\gamma}_{B}$, where $B\in\mathfrak{B}(\A)$. Given an irreducible interval $[\s_1,\s_2]$ in $\silt(\A)$, let $s_1$ and $s_2$ be the silting objects in $H_0\A$ determining $\s_1$ and $\s_2$, respectively, and let $r$ be the rigid object such that $H_0(\s_1\cap \s_2)=\add(r)$ (see \Cref{AT.5.4}).
Consider the following diagram in $\tcmc{\A}$:
\begin{equation*}
		\begin{tikzcd}[column sep=4.75em]
			\A \arrow[r,"\s_1\cap \s_2"] & \A/(\s_1\cap \s_2) \arrow[r,yshift=0.25em,"\varpi_{\s_1\cap \s_2}(\s_1)"] \arrow[r,yshift=-0.25em,"\varpi_{\s_1\cap \s_2}(\s_2)"'] & \mathscr{O}.
		\end{tikzcd}
		\end{equation*}
The two composites are the silting subcategories $\s_1$ and $\s_2$. Since $s_1$ and $s_2$ are the only rigid objects in $H_0\A$ properly containing $r$ as a direct summand, the morphisms $\varpi_{\s_1\cap \s_2}(\s_1)$ and $\varpi_{\s_1\cap \s_2}(\s_2)$ are the only ones in $\tcmc{\A}$ with domain $\A/(\s_1\cap \s_2)$. Consequently, we have that $\A/(\s_1\cap \s_2)\in \mathfrak{B}(\A)$. 
The construction of the (elementary) homotopy 
\begin{equation*}
		\begin{tikzcd}[column sep=6em,row sep=3em]\label{eq:picgroup_pres2.5.two-hom}
\mathscr{O} \arrow[d, equal] & \A \arrow[l, "\s_2"'] \arrow[r, "\s_1"] \arrow[d, "\s_1\cap \s_2"']                                   & \mathscr{O} \arrow[d, equal] \\
\mathscr{O}                    & \A/(\s_1\cap \s_2) \arrow[l, "\varpi_{\s_1\cap\s_2}(\s_2)"'] \arrow[r, "\varpi_{\s_1\cap\s_2}(\s_1)"] & \mathscr{O}          
\end{tikzcd}
	\end{equation*}
concludes the proof.
\end{proof}

\section{Future work}\label{subsec:future}

\subsection{} Let $\Gamma$ be a group. A functor of the form $\c \to \Gamma$ is called a \textit{group functor}, where we regard $\Gamma$ as a groupoid with a single object. If $\A$ is a 0-Auslander exact dg $\k$-category for which the picture category $\tcmc{\A}$ of $\A$ is defined,
we say that a $\A$ \textit{admits a faithful group functor} if there exists a group $\Gamma$ and a group functor $\tcmc{\A} \to \Gamma$ which is faithful.

\subsection{}\label{Igu14.3.4} The existence of faithful group functors, along with the cubical structure of picture categories (see \Cref{prop:HI18.2.14}), can be useful to determine the homotopy type of picture spaces. Indeed, a cubical category $\c$ is locally $\mathrm{CAT}(0)$ precisely when the conditions \hyperref[item:I1]{\textbf{(I1)}}--\hyperref[item:I3]{\textbf{(I3)}} below are met \cite[Propositions 3.4 and 3.7]{Igu14}. In particular, these conditions ensure that the classifying space $\classspace\c$ is an Eilenberg--MacLane space of type $\K(\pi,1)$ \cite{Gro87}, which is to say that the homotopy groups $\pi_i(\classspace\c,\ast)$ are trivial for all $i\geq 2$. 
\begin{enumerate}
	\item[\textbf{(I1)}]\label{item:I1} A set of $\ell$ rank 1 morphisms $\{X \to Y_j \}$ in $\c$ forms the set of first factors of a rank $\ell$ morphism in $\c$ if and only if each subset of cardinality 2 forms the set of first factors of a rank 2 morphism in $\c$.
	\item[\textbf{(I2)}]\label{item:I2} A set of $\ell$ rank 1 morphisms $\{Y_j \to X \}$ in $\c$ forms the set of last factors of a rank $\ell$ morphism in $\c$ if and only if each subset of cardinality 2 forms the set of last factors of a rank 2 morphism in $\c$.
	\item[\textbf{(I3)}]\label{item:I3} There exists a faithful group functor $\c\to\Gamma$, for some group $\Gamma$.
\end{enumerate}
Since the picture space $\picspace{\A}$ is connected, as was argued in \Cref{tcmc_str}, we can thus determine its homotopy type if the conditions \hyperref[item:I1]{\textbf{(I1)}}--\hyperref[item:I3]{\textbf{(I3)}} above are met.

\subsection{}
It is straightforward to show that the condition \hyperref[item:I1]{\textbf{(I1)}} above holds for the picture category $\tcmc{\A}$ of a 0-Auslander exact dg $\k$-category $\A$ with a silting object. Less trivially, we will prove that \hyperref[item:I3]{\textbf{(I3)}} holds for $\tcmc{\A}$ when $\k$ is the field of complex numbers \cite{BorMot}. This will be achieved with the aid of motivic Hall algebras.
We will follow Bridgeland \cite{Bri17} when defining the completed motivic Hall algebra $\widehat{\mathrm{H}}(\A)$ of $\A$, as well as the pro-unipotent group $\widehat{G}_{\mathrm{Hall}}(\A)$ corresponding to the pro-nilpotent Lie $\mathbb{Q}(t)$-algebra $\widehat{\mathrm{H}}_{>0}(\A)$ inside $\widehat{\mathrm{H}}(\A)$. We provide a faithful group functor $\tcmc{\A} \to \widehat{G}_{\mathrm{Hall}}(\A)$, whereby proving a motivic counterpart of a conjecture of Hanson--Igusa \cite[Conjecture 5.10]{HI21p}. 

\bibliographystyle{IEEEtranSA}
\bibliography{0Aus}

\end{document}